\documentclass[10pt]{amsart}
\usepackage{amscd, amssymb}
\newtheorem{Thm}{Theorem}[subsection]
\newtheorem{Conj}[Thm]{Conjecture}
\newtheorem{Prop}[Thm]{Proposition}
\newtheorem{Def}[Thm]{Definition}
\newtheorem{Def/Thm}[Thm]{Definition/Theorem}
\newtheorem{Cor}[Thm]{Corollary}
\newtheorem{Lemma}[Thm]{Lemma}

\theoremstyle{remark}
\newtheorem{Rmk}[Thm]{Remark}

\numberwithin{equation}{subsection}
\newcommand{\ti }{\times}
\newcommand{\ot }{\otimes}
\newcommand{\ra }{\rightarrow}

\newcommand{\lra }{\longrightarrow}
\newcommand{\Hom }{{\mathrm{Hom}}}

\newcommand{\Pic}{{\mathrm{Pic}}}

\newcommand{\cO}{{\mathcal{O}}}
\newcommand{\cM}{{\mathcal{M}}}
\newcommand{\cL}{{\mathcal{L}}}
\newcommand{\cE}{{\mathcal{E}}}
\newcommand{\cV}{{\mathcal{V}}}
\newcommand{\fE}{{\mathfrak{E}}}
\newcommand{\G}{{\bf G}}
\newcommand{\T}{{\bf T}}
\newcommand{\bS}{{\bf S}}
\newcommand{\bW}{{\bf W}}
\newcommand{\PP }{{\mathbb P}}
\newcommand{\QQ }{{\mathbb Q}}
\newcommand{\CC }{{\mathbb C}}
\newcommand{\ZZ }{{\mathbb Z}}

\newcommand{\bt }{{ {t}}}
\newcommand{\ka }{{\alpha}}

\let\wt\widetilde
\newcommand{\onabla}{{}^\omega\!\nabla}
\DeclareMathOperator{\sgn}{sgn}

\begin{document}
\title{The Abelian/Nonabelian Correspondence and Frobenius Manifolds}

\begin{abstract}
We propose an approach via Frobenius manifolds to the study (began
in \cite{BCK2}) of the relation between rational Gromov-Witten
invariants of nonabelian quotients $X//\G$ and those of the
corresponding ``abelianized'' quotients $X//\T$, for $\T$ a
maximal torus in $\G$. The ensuing conjecture expresses the
Gromov-Witten potential of $X//\G$ in terms of the potential of
$X//\T$. We prove this conjecture when the nonabelian quotients
are partial flag manifolds.
\end{abstract}

\author{Ionu\c t Ciocan-Fontanine}
\address{School of Mathematics, University of Minnesota,
Minneapolis MN, 55455, USA}
\email{ciocan@math.umn.edu}

\author{Bumsig Kim}
\address{School of Mathematics, Korea Institute for Advanced Study,
207-43 Cheongnyangni 2-dong, Dongdaemun-gu, Seoul, 130-722, Korea}
\email{bumsig@kias.re.kr}

\author{Claude Sabbah}
\address{UMR 7640 du C.N.R.S., Centre de math\'ematiques Laurent Schwartz,
\'Ecole polytechnique, F-91128 Palaiseau cedex, France}
\email{sabbah@math.polytechnique.fr}

\maketitle

\section{Introduction}
\subsection{} The paper \cite{BCK2}
conjectures a correspondence between the genus zero
Gro\-mov-Witten invariants of nonsingular projective GIT quotients
$X//\G$ and $X//\T$, for $\G$ a complex reductive Lie group with a
linearized action on a projective manifold $X$ and $\T$ a maximal
torus in $\G$. The correspondence expresses (descendant)
Gromov-Witten invariants of $X//\G$ in terms of Gromov-Witten
invariants of $X//\T$ {\em twisted} by (the top Chern class of) a
certain decomposable vector bundle on $X//\T$.

Our main goal in this paper is to give a natural reformulation of
the correspondence in terms of the Frobenius structures describing
the (big) quantum cohomology rings $QH^*(X//\G,\CC)$ and
$QH^*(X//\T,\CC)$. This is accomplished in section \ref{Frobenius
correspondence}. To explain it, recall that a given cohomology
class $\sigma\in H^*(X//\G)$ can be lifted to a class
$\tilde{\sigma}$ (of the same degree) in the Weyl group invariant
subspace $H^*(X//\T)^{\bW}$. Such a lift is not unique, however,
if $\omega$ is the fundamental $\bW$-anti-invariant class, then
$\tilde{\sigma}\cup\omega$ is uniquely determined by $\sigma$.
Moreover, by results of Ellingsrud and Str\o mme when $X=\PP^N$,
and later Martin in full generality, this identification respects
cup products:
$$(\widetilde{\sigma\cup_{X//\G}\sigma'})\cup\omega
=\tilde{\sigma}\cup(\tilde{\sigma}'\cup\omega) \in H^*(X//\T).$$
A naive guess might be that the identification respects quantum products as well, that is,
$$(\widetilde{\sigma\star_{X//\G}\sigma'})\cup\omega
=\tilde{\sigma}\star_{X//\T}(\tilde{\sigma}'\cup\omega),$$ after
an appropriate specialization of quantum parameters. Indeed, as
shown in \cite{BCK1}, Theorem 2.5, this is the case for {\it
small} quantum products when $X//\G$ is a Grassmannian. At the
level of Gromov-Witten invariants, this would translate into an
appealing identity of the form
\begin{equation}\label{naive}
\langle\sigma_1,\sigma_2\dots,\sigma_{n-1},\sigma_n\rangle_{0,n,\beta}^{X//\G}
=\pm\sum_{\tilde{\beta}\mapsto\beta}
\langle\tilde{\sigma}_1,\dots,\tilde{\sigma}_{n-2},\tilde{\sigma}_{n-1}\cup\omega,
\tilde{\sigma}_n\cup\omega\rangle_{0,n,\tilde{\beta}}^{X//\T}.
\end{equation}
It is not hard to convince oneself, however, that this fails for
big quantum cohomology (already for the Grassmannian
$Grass(2,4)$), and that it has no reason to be true in general
even for small quantum cohomology. Instead, we conjecture a
generalization to quantum cohomology as follows:

\noindent {\em Fix}
a lifting $\;\widetilde{\bullet}\;$ of $H^*(X//\G)$ to a subspace
$U\subset H^*(X//\T)^{\bW}$. Let $\{t_i\}$ be the coordinates on
$H^*(X//\G)$, corresponding to a choice of basis, and let
$\{\tilde{t}_i\}$ be the coordinates on $U$ corresponding to the
lifted basis. Let $N(X//\G)$ and $N(X//\T)$ be the Novikov rings
for the two quotients.

The quantum product in $QH^*(X//\G,\CC)$ is
a $N(X//\G)[[t]]$-linear product on $H^*(X//\G,\CC)\otimes_\CC
N(X//\G)[[t]]$, while the quantum product in $QH^*(X//\T,\CC)$ is
a $N(X//\T)[[\tilde{t},y]]$-linear product on $H^*(X//\T,\CC)\otimes_\CC
N(X//\T)[[\tilde{t},y]]$, where $(\tilde{t},y)$ is an extension of
$\tilde{t}$ to coordinates on the entire $H^*(X//\T,\CC)$.

There is a natural specialization of Novikov variables
$p:N(X//\T)\ra N(X//\G)$ which takes into account that there are
more curve classes on $X//\T$. We denote by $\lq\lq{\star}"$ the
quantum product on ${X//\T}$ with the Novikov variables
specialized via $p$. Given $\sigma,\sigma'\in H^*(X//\G)$, there
are classes $\xi,\xi'\in U\otimes _{\CC}
N(X//\G)[[\;\tilde{t}\;]]$, uniquely determined by
$\xi{\star}\omega=\tilde{\sigma}\cup\omega$ and
$\xi'{\star}\omega=\tilde{\sigma}'\cup\omega$ respectively.

\bigskip
\noindent{\bf Conjecture.} {\it There is an equality
$$\left((\widetilde{\sigma\star_{ X//\G}\sigma'})\cup\omega\right) (t)=
\xi\star \xi'\star \omega \; (\tilde{t},0)=\left(\xi\star
(\tilde{\sigma}'\cup\omega)\right)(\tilde{t},0),$$ after an
explicit change of variable $\tilde{t}=\tilde{t}(t)$. }

\bigskip

At the level of Gromov-Witten invariants, the Conjecture says that
the right-hand side of the naive formula (\ref{naive}) receives a
correction term which is a sum of products of invariants of $X//\T$
of the same type (see the Appendix for a discussion and some examples).

We should warn the reader that the above formulation is a
translation of the actual Conjecture \ref{conj} in the body of the
paper, which is stated in the conceptual framework of Frobenius
structures. It is in this framework that one is naturally lead to
the conjecture. Indeed, if $N$ is the formal germ of the affine
space over $N(X//\G)$ associated to the subspace $U$, the general
machinery of the infinitesimal period mapping in the theory of
Frobenius-Saito structures (see e.g., \cite{Sabbah}) gives a
canonical Frobenius manifold structure on $N$. It is induced by
the primitive homogeneous section $\omega$ of the (trivial) bundle
with fiber the anti-invariant subspace $H^*(X//\T)^a$ over $N$,
together with the restriction to this bundle (in an appropriate
sense) of the Frobenius structure on $H^*(X//\T)$. Our conjecture
says that this new Frobenius manifold is identified with the
Frobenius manifold given by the Gromov-Witten theory of $X//\G$.
The new flat metric $^\omega g$ on the sheaf $\Theta_N$ of vector
fields
 satisfies $$^\omega g(\tilde{\sigma},\tilde{\sigma}')=
g(\tilde{\sigma}\star\omega,\tilde{\sigma}'\star\omega).$$ It
follows that the coordinates $\{\tilde{t}_i\}$ on $N$ provided by
lifting are {\it not} flat for the new Frobenius structure, or,
equivalently, the liftings $\tilde{\sigma}$ are not horizontal
vector fields. The vector fields $\xi$, $\xi'$ appearing in the
statement of the conjecture are precisely the horizontal vector
fields corresponding to $\sigma$, $\sigma'$ under the
identification of flat coordinates of Frobenius structures. This
identification of coordinates is the change variable
$\tilde{t}=\tilde{t}(t)$.

In fact, we treat a more general situation in \S \ref{Frobenius
correspondence}, by considering the {\em equivariant}
Gromov-Witten theories in the presence of compatible actions of an
additional torus $\bS$ on $X//\G$ and $X//\T$. The corresponding
Frobenius structures are more general than the ones considered in
\cite{Sabbah}, as they lack an Euler vector field. However, a
suitable modification of the notion of Euler vector field allows
the application of the theory of infinitesimal period mappings in
this case as well. We give an exposition of the relevant facts in
\S \ref{equiv Frobenius}--\ref{finite group}.

This generalization is needed in \S \ref{the proof}, where we prove,
by using reconstruction theorems for Gromov-Witten
invariants (extended to the equivariant setting), that the conjecture above can
be reduced in many cases to the abelian/nonabelian correspondence for {\it small}
$J$-functions from \cite{BCK2}. In particular, the following result is obtained:

\bigskip

\noindent{\bf Theorem.}\;\;{\it Let $Fl=Fl(k_1,\dots,k_r,n)$ be
the flag manifold parameterizing flags of subspaces
$\{\CC^{k_1}\subset\dots\subset\CC^{k_r}\subset \CC^n\}$, viewed
as a GIT quotient $\PP^l//\G$ for appropriate $l,\G$. Denote by
$Y$ the toric variety which is the corresponding abelian quotient
$\PP^l//\T$ (cf. \cite{BCK2}). Then the conjecture is true for the
pair $(Fl,Y)$. }

\bigskip

The Theorem implies that the genus zero Gromov-Witten invariants
of a flag manifold (with any number of insertions) can be
expressed in terms of Gromov-Witten invariants of the associated
toric variety $Y$. In an appendix we write down explicit formulae
in the simplest case of the Grassmannian $Grass(k,n)$, for which
the abelian quotient is the product of $k$ copies of $\PP^{n-1}$.

In section \ref{J-function}, we obtain an equivalent formulation
(\ref{J-I correspondence})
 of the conjecture in terms of
(big) $J$-functions of $X//\T$ and $X//\G$. It generalizes
Conjecture 4.3 of \cite{BCK2} and, by the above theorem, it holds
for type $A$ flag manifolds.

Finally, in section \ref{others} we extend the abelian/nonabelian
correspondence to include Gromov-Witten invariants
with an additional twist by homogeneous
vector bundles. As an application, we
describe the $J$-function of a generalized flag manifold for a simple
complex Lie group of type $B$, $C$, or $D$ as the twisted
$J$-function of the abelianization of the corresponding flag
manifold of type $A$.
\medskip

\subsection{Acknowledgements} Ciocan-Fontanine's work was partially supported
by the NSF grant DMS-0303614. Part of the final writing of the
paper was done during a visit by Ciocan-Fontanine to KIAS, whose
support and hospitality are gratefully acknowledged. Kim's work is
supported by KOSEF grant R01-2004-000-10870-0. Kim thanks staffs
at \'Ecole Polytechnique for their warm hospitality during his
visit. Sabbah thanks KIAS for providing him with excellent working
conditions during his visit.

\section{Preliminaries on Frobenius structures}\label{preliminaries}

\subsection{Formal Frobenius manifolds from Gromov-Witten theory}
Let $R$ be a $\CC$-algebra and let $K$ be a free $R$-module of rank $m$.
We think of $K$ as the affine $m$-space over $R$ (precisely, the spectrum
of the symmetric algebra of the dual module). Let
$M:=Spf(R[[K^{\vee}]])$ be the formal completion of $K$ at the origin.
$M$ is a formal manifold over $R$. We denote by $\Theta_M$ its formal relative
tangent sheaf over $R$. Note that it
is canonically identified with $K\otimes_R\cO_M$.

\begin{Def}\label{Frob} The data $(M, \star, g, e ,\fE)$ is
called a (conformal, even) formal Frobenius manifold over $R$
if the following properties hold:

\begin{itemize}

\item $g$ is an $R$ - linear, nondegenerate pairing such
that its metric connection $\nabla$ is flat

\item $\star$ is an $R$ - linear, associative, commutative
product on $\Theta_M$

\item $e$ is a formal vector field on $M$ over $R$
which is the identity for the product $\star$, and such that
$\nabla e =0$

\item $\nabla c$ is symmetric, where the tensor $c$ is defined by
$c(u,v,w)=g(u\star v,w)$

\item $\fE$ is a formal vector field on $M$ over $R$
satisfying
$$\cL_{\fE}(g)=Dg,\;\;\; \cL_{\fE}(\star)=\star,\;\;\; \cL_{\fE}(e)=-e,$$
where $\cL_{\fE}$ denotes the Lie derivative and $D\in \CC$ is a constant.

\end{itemize}
\end{Def}

The fourth condition implies that there is a formal function $F$
on $M$ (the {\em potential} of the Frobenius manifold) such that
the tensor $c$ is given by the third derivatives of $F$ in flat
coordinates, and then associativity of $\star$ translates into the
WDVV equations for $F$. The vector field $\fE$ is called an {\em
Euler vector field}.

We recall here the formal Frobenius manifold structures determined by the
genus zero GW-theories (ordinary and equivariant) of a projective manifold endowed with an action of an
algebraic complex torus ${\bf S}\cong ({\CC}^*)^\ell$.
Detailed expositions can be found in \cite{LP2} or \cite{Manin}, to which
we refer the reader.

Let $Y$ be a smooth projective variety over $\CC$. We assume for
simplicity that $H_2(Y,\ZZ)$ is torsion-free and that the odd
cohomology $H^{2*+1}(Y,\CC)$ vanishes. We denote by $N(Y)$ the
{\em Novikov ring} of $Y$. It can be described as the
$\CC$-algebra of \lq\lq power series" $\{\sum_{\beta\in NE_1}
c_\beta Q^\beta | c_\beta\in\CC\}$, where $NE_1\subset H_2(Y,\ZZ)$
is the semigroup of effective curve classes.

The genus zero Gromov-Witten theory of $Y$ determines a formal Frobenius manifold
over $R=N(Y)$. We take
$$K=N(Y)\otimes_{\CC} H^*(Y,\CC),$$
so that
$M= Spf (N(Y)[[K^{\vee}]])$.
%and the global vector fields $\Gamma(M,\Theta_M)$
%are identified canonically with $K$.
The metric $g$ is given by the intersection pairing:
$$g(\gamma,\gamma')=\int_Y\gamma\cup\gamma'.$$
Let $\{ 1=\gamma_0,\gamma_1,\dots,\gamma_r,\gamma_{r+1},
\dots, \gamma_{m-1}\}$
be a basis of $H^*(Y,\CC)$ consisting of integral homogeneous
classes, such that $\gamma_1,\dots,\gamma_r$ form
a basis of $H^2$.
We write
$\sigma=\sum t_i\gamma_i$ for a general cohomology class on $Y$.
The functions $t_i$ give flat coordinates on $M$.
A potential function is defined using the genus zero Gromov-Witten invariants of $Y$
$$F(Q,t):=\sum_{\beta\in NE_1}\sum_{n\geq 0}Q^\beta\frac{1}{n!}
\langle\underbrace{\sigma,\dots, \sigma}_{n} \rangle_{0,n,\beta},$$
where the unstable terms with $\beta=0,\;n\leq 2$ are omitted in the sum.
The tensor $c$ is given in flat coordinates by
$$c_{ijk}=\partial_{t_i}\partial_{t_j}\partial_{t_k} F$$ and the product $\star$ is called the
{\em big quantum product}.
The unit vector field $e$ is given by the class $\gamma_0=1$.

The following notation is customary:
$$\langle\langle \sigma_1,\dots ,\sigma_r\rangle\rangle=
\sum_{\beta\in NE_1}\sum_nQ^\beta\frac{1}{n!}\langle\sigma_1,\dots ,\sigma_r,
\underbrace{\sigma,\dots \sigma}_{n} \rangle_{0,n+r,\beta},$$
where $\sigma_j\in H^*(Y,\CC)$ are given cohomology classes and
$\sigma=\sum t_i\gamma_i$ is the general element in $H^*(Y,\CC)$ (so that $\langle\langle\;\;\rangle\rangle=F$).
We extend this double bracket $\cO_M$-linearly to general vector fields
$\sigma_1,\dots ,\sigma_r$.
It is easy to see that for any vector field $\xi$ we have
$$\nabla_\xi(F)= \langle\langle\xi\rangle\rangle.$$
In particular, since $\nabla_{\partial_{t_i}}\partial_{t_j}=0$, the quantum product
can be written in our chosen basis
$$\gamma_i\star\gamma_j=\sum_k\langle\langle\gamma_i,\gamma_j,\gamma_k\rangle\rangle
\gamma_k^{\vee}$$ where $\gamma_k^{\vee}=\sum_lg^{kl}\gamma_l$ with
$(g^{kl})$ the inverse matrix of the metric $g$.

The divisor axiom for Gromov-Witten invariants implies that the Gromov-Witten potential has
the special form
\begin{eqnarray} \label{qc-type} F&=&F_{cl}+\sum_{\beta\in NE_1,\beta\neq 0}
Q^{\beta}e^{\beta\cdot t_{\mathrm{small}}}F_{\beta},
\end{eqnarray}
with $F_{cl}$ a cubic polynomial in the $t_i$'s and $F_{\beta}\in\CC[[t_{r+1},\dots,t_{m-1}]]$
formal power series in the
{\em non-divisorial}
coordinates. Here we use the notation $\beta\cdot t_{\mathrm{small}}$ for the intersection
index of $\beta$ with the general $H^2$-class,
$$\beta\cdot t_{\mathrm{small}} :=\int_{\beta}\sum_{i=1}^rt_i\gamma_i.$$
We will also use the notation $F_q$ for $F-F_{cl}$.

Assume now that $Y$ is acted upon by the torus ${\bf S}\cong ({\CC}^*)^\ell$.
The equivariant cohomology $H^*_{\bf S}(Y,\CC)$ is a module over the
polynomial ring
$$H^*_{\bf S}(pt)=H^*(B{\bf S})\cong\CC [\lambda _1,...,\lambda _\ell],$$
and it is in fact a free module by \cite{Ginzburg}. Taking $$R=N(Y)[\lambda]:=N(Y)\otimes_{\CC}
\CC [\lambda _1,...,\lambda _\ell]$$ and
$$K_\bS=R\otimes_{\CC[\lambda _1,...,\lambda _\ell]} H^*_{\bf S}(Y,\CC)$$
we get similarly a formal Frobenius manifold over $R$. The metric $g$ is now
given by the ($\CC [\lambda _1,...,\lambda _\ell]$-valued)
equivariant intersection pairing,
while in $F$ the GW-invariants are replaced by their ${\bf S}$-equivariant
counterparts. The unit vector field and equivariant big quantum product are
obtained analogously.

Localization with respect to ${\bf S}$
determines yet another Frobenius structure. Consider the
localization of $H^*(B{\bf S})$, i.e., the field of fractions
$\CC(\lambda_1,\dots, \lambda_{\ell})$, and set
$$N(Y)[\lambda]_{(\lambda)}=N(Y)\otimes_{\CC}\CC(\lambda_1,\dots, \lambda_{\ell})$$
$$K_\bS^*=N(Y)[\lambda]_{(\lambda)}\otimes_{N(Y)} K_{\bf S}.$$
Taking $M=Spf(N(Y)[\lambda]_{(\lambda)}[[K_\bS^{*\vee}]])$
with the localized equivariant metric, potential function,
and unit vector field determines a formal Frobenius manifold over
$N(Y)[\lambda]_{(\lambda)}$(in other words, we simply consider the Frobenius
structure induced by base change via $N(Y)\ra N(Y)[\lambda]_{(\lambda)}$).

In both the equivariant and localized equivariant cases the potential function in flat
coordinates $t$ has the special form (\ref{qc-type}), with $F_{\beta}\in
\CC[\lambda][[t_{r+1},\dots,t_{m-1}]]$.

Finally, we discuss the Euler vector fields.
The Frobenius manifold defined by the (nonequivariant) Gromov-Witten theory
of $Y$ is conformal: the Euler vector field (with $D=2-\dim (Y)$) is explicitly
$$\fE=\sum _{i=0}^{m-1} (1-\frac{\mathrm{ cdeg} \gamma
_i}{2}){t}_i
\partial _{{t}_i}+c_1(TY).
$$
Here \lq\lq cdeg" is the cohomological degree.

Consider the ${\bf S}$-equivariant version
of this vector field

$$\fE=\sum _{i=0}^{m-1} (1-\frac{\mathrm{ cdeg} \gamma
_i}{2}){t}_i
\partial _{{t}_i}+c_1^{\bf S}(TY)
$$
with $\gamma_i$'s now forming a basis of $H^*_{\bf S}(Y,\CC)$ over $H^*(B{\bf S})$.
$\fE$ does not
give a conformal structure (because equivariant Gromov-Witten invariants
do not satisfy a dimension constraint). Nevertheless, we consider
a variant of the Euler vector field in this context as well, by relaxing the
requirement of linearity over $N(Y)\otimes \CC[\lambda]$ and
will define below an Euler vector field $\mathfrak{E}_\bS$ as an
$N(Y)$-derivation on $\cO_M$ (that is, an $N(Y)$-derivation of $K_\bS^\vee$ into itself).
%$N(Y)[\lambda]\otimes _{\CC [\lambda]} H^*_\bS(Y)^\vee$ into itself.
%$N(Y)[\lambda]\otimes _{\CC [\lambda]} H^*_\bS(Y)^\vee$ (See
%Hartshorne for the terminology)). First choose a cohomologically
%homogeneous basis $\{\gamma _i\}$ of $H^*_T(Y)$ over $\CC
%[\lambda]$ and let
The flat coordinates $\{t_i\}$, together with $\{\lambda_1,\dots,\lambda_{\ell}\}$
form a coordinate system on $M$ over $Spec(N(Y))$. Therefore
%(or, more precisely, $1\otimes t_i$)
%are then the
%corresponding dual elements of $N(Y)[\lambda]\otimes _{\CC
%[\lambda]} H^*_\bS(Y)^\vee$ (where $H^*_T(Y)^\vee =
%\mathrm{Hom}_{\CC [\lambda ]} (H^*_T(Y),\CC [\lambda ])$. Now
%using the isomorphism $N(Y)[\lambda][[t_0,...,t_{m-1}]]
%\cong N(Y)[\lambda]\otimes H^*_T(Y)^\vee$
%we define the relative tangent vector field (that is,
%$N(Y)[\lambda]$-derivation) $\sum _i (1-\frac{\mathrm{cdeg} \gamma
%_i}{2})t_i\partial _{t_i} + c_1^{\mathbf{S}}(TY)$ and define a
%$N(Y)$-derivation $\cE _{\bf S} :=\sum _{i=1}^l \lambda
%_i\partial_{\lambda _i}$ on $M$. Then set $\mathfrak{E}_{\bf S}:=
%\mathfrak{E}+\cE _S$ as a $N(Y)$-derivation on $M$. Note that the
%definition of the Euler vector field does not depend on the
%choices of the homogeneous bases $\{\gamma _i\}$ since
%$\mathfrak{E}_{\bf S} (s) = (\deg s) s$ if $s$ is a homogeneous
%function on $M$ over $N(Y)$.
$$\cE_{\bf S}:=\sum_{i=1}^{\ell}\lambda_i\partial_{\lambda_i}$$
is a well-defined \lq\lq absolute" vector field (i.e.,
$N(Y)$-linear derivation) and acts by Lie bracket on the relative
vector fields $\Theta_M$. Put
$$\fE_{\bf S}:=\fE+\cE_{\bf S}.$$
If $\eta\in \Theta_M$
is any relative vector field, then the commutator $[\fE_{\bS},\eta]$ is also in $\Theta_M$.
Hence Lie derivatives of tensors on $\Theta_M$ are well defined.
The vector field $\fE_{\bf S}$ will still satisfy the conditions in Definition
\ref{Frob}, again with $D=2-\dim (Y)$. The same $\fE_{\bf S}$ will be used for the
localized structure as well.

We have
$$\cL_{\fE_\bS}(\partial_{t_i})=(-1+\frac{{\mathrm{cdeg}}\gamma_i}{2})\partial_{t_i},\;\;\;
\cL_{\fE_\bS}(\lambda _i)= \lambda _i ,\;\;\;
\cL_{\fE_\bS}(F)=(3-\dim(Y))F.$$

\subsection{$\bS$-Equivariant Frobenius manifolds over $R$}\label{equiv
Frobenius} We extend the construction of Frobenius manifold
through an infinitesimal period mapping to the previous setting.
Let $M$ be as above, with $\cO_M=R[[K^{*\vee}_\bS]]$ and
$R=N(Y)[\lambda]$ or $R=N[Y][\lambda]_{(\lambda)}$. Let $E$ be a
free $\cO_M$-module of finite rank. An \emph{$\bS$-conformal
connection} on $E$ consists of a pair
$\wt\nabla=(\nabla,\wt\nabla_{\cE_\bS})$, where $\nabla$ is an
$R$-connection on $E$ and $\wt\nabla_{\cE_\bS}$ is a $N(Y)$-linear
derivation satisfying $\wt\nabla_{\cE_\bS}(\varphi
e)=\cL_{\cE_\bS}(\varphi)e+\varphi\wt\nabla_{\cE_\bS}e$ for any
$e\in E$ and $\varphi\in\cO_M$. We say that $\wt\nabla$ is
\emph{flat} if $\nabla$ is flat and for any vector field
$\xi\in\nobreak\Theta_M$,
$[\wt\nabla_{\cE_\bS},\nabla_\xi]=\nabla_{[\cE_\bS,\xi]}$. In
coordinates $(t_i)$ defined from an $N(Y)$-basis of $K$, the
previous condition is equivalent to the pairwise commutation of
the operators $\nabla_{\partial_{t_i}}$ and $\wt\nabla_{\cE_\bS}$.
Such a connection $\wt\nabla$ extends in a natural way to a
similar object on $\hom_{\cO_M}(E,E)$.

An $\bS$-equivariant \emph{pre-Saito structure}
$(M,E,\wt\nabla,\Phi,R_0,g)$ of weight $w$ over $M$ consists of
\begin{itemize}
\item
a free $\cO_M$-module $E$ of finite rank with a flat $\bS$-conformal connection
$\wt\nabla$,
\item
$\cO_M$-linear morphisms $\Phi:\Theta_M\otimes_{\cO_M}E\to E$ and $R_0:E\to E$,
\item
an $\cO_M$-bilinear form $g$ on $E$,
\end{itemize}
satisfying, when expressed in coordinates $(t_i)$, the following relations for
all $i,j$:
\begin{gather*}
\nabla_{\partial_{t_i}}(\Phi_{\partial_{t_j}})=\nabla_{\partial_{t_j}}(\Phi_{\partial_{t_i}}),\quad
[\Phi_{\partial_{t_i}},\Phi_{\partial_{t_j}}]=0,\quad
[R_0,\Phi_{\partial_{t_i}}]=0,\\
\Phi_{\partial_{t_i}}-\wt\nabla_{\cE_\bS}(\Phi_{\partial_{t_i}})+\nabla_{\partial_{t_i}}(R_0)=0,\\
\nabla(g)=0,\quad
\wt\nabla_{\cE_\bS}(g)=-wg,\quad\Phi_{\partial_{t_i}}^*=\Phi_{\partial_{t_i}},\quad
R_0^*=R_0,
\end{gather*}
where ${}^*$ means taking the $g$-adjoint and $\wt\nabla_{\cE_\bS}(g)$ is
defined as usual by the formula
$\wt\nabla_{\cE_\bS}(g)(\xi,\eta)=\cL_{\cE_\bS}\big(g(\xi,\eta)\big)
-g(\wt\nabla_{\cE_\bS}\xi,\eta)-g(\xi,\wt\nabla_{\cE_\bS}\eta)$.

The pull-back of an $\bS$-equivariant pre-Saito structure by a
morphism $f:N\to M$ is well-defined only for morphisms $f^*$ which
commute with $\cL_{\cE_\bS}$.

The definition of an \emph{$\bS$-equivariant Frobenius manifold
over $R$} is a variant of Definition \ref{Frob}: With the same
data $(M,\star,g,e,\fE)$, we set $\fE_\bS=\fE+\cE_\bS$, and we
replace the homogeneity conditions by the following ones:
\[
\cL_{\fE_\bS}(g)=Dg,\quad \cL_{\fE_\bS}(\star)=\star,\quad\cL_{\fE_\bS}(e)=-e.
\]

Let $(M,E,\wt\nabla,\Phi,R_0,g)$ be an $\bS$-equivariant
\emph{pre-Saito structure} of weight $w$ and let $\omega$ be a
$\nabla$-horizontal section of $E$. It defines an $\cO_M$-linear
morphism $\varphi_\omega:\Theta_M\to E$ by
$\xi\mapsto-\Phi_\xi(\omega)$. We say that such a section $\omega$
of $E$ is
\begin{enumerate}
\item
\emph{primitive} if the associated period mapping $\varphi_\omega:\Theta_M\to
E$ is an isomorphism,
\item
\emph{homogeneous} of degree $q\in\CC$ if $\wt\nabla_{\cE_\bS}\omega=q\omega$.
\end{enumerate}

The data of an $\bS$-equivariant pre-Saito structure and of a
homogeneous primitive section $\omega$ is called an
$\bS$-equivariant Saito structure. As in \cite[\S4.3]{Sabbah} and
following K.~Saito, we obtain the following results.

If $\omega$ is primitive and homogeneous, $\varphi_\omega$ induces a flat
torsionless $R$-connection $\onabla:=\varphi_\omega^{-1}\nabla\varphi_\omega$
on $\Theta_M$, and an associative and commutative $\cO_M$-bilinear
product~$\star$ by $\xi\star\eta=-\Phi_\xi(\varphi_\omega(\eta))$, with
$e=\varphi_\omega^{-1}(\omega)$ as unit, and $\onabla e=0$. Moreover, $\onabla$
is the metric connection attached to the metric ${}^\omega\!g$ on $\Theta_M$
obtained from $g$ through $\varphi_\omega$, and $\onabla$ is $\bS$-conformal
and flat as such, setting $\wt{\onabla}_{\cE_\bS}=\varphi_\omega^{-1}\circ
\wt\nabla_{\cE_\bS}\circ\varphi_\omega-\mathrm{Id}$.

The Euler field is $\fE=\varphi_{\omega}^{-1}(R_0(\omega))$. It is therefore a
section of $\Theta_M$. We have
$\onabla\fE=\cL_{\cE_\bS}-\wt\onabla_{\cE_\bS}+q\mathrm{Id}$. In particular,
$\onabla(\onabla\fE)=0$.

If we put $D=2q+2-w$, and if we set as above $\fE_\bS=\fE+\cE_\bS$, we get
\[
\cL_{\fE_\bS}(e)=-e,\quad
\cL_{\fE_\bS}(\star)=\star,\quad\cL_{\fE_\bS}({}^\omega\!g)=D\cdot{}^\omega\!g.
\]

Given an $\bS$-equivariant pre-Saito structure
$(M,E,\wt\nabla,\Phi,R_0,g)$ of weight~$w$, the datum of a
homogeneous primitive section $\omega$ of $E$ having weight $q$
induces on $M$, through $\varphi_\omega$, the structure o f a
$\bS$-equivariant Frobenius manifold of weight $D=2q+2-w$.

Conversely, any $\bS$-equivariant Frobenius manifold
$(M,\star,g,e,\fE)$ defines an $\bS$-equivariant pre-Saito
structure $(M,\Theta_M,\wt\nabla,\Phi,R_0,g)$ having $e$ as
homogeneous primitive form.

For instance, to give the correspondence
$(M,\star,g,e,\fE)\mapsto(M,\Theta_M,\wt\nabla,\Phi,R_0,g)$ we take $\nabla$
to be the
Levi-Civita connection of $g$, and
\[
\wt\nabla_{\cE_\bS}=\mathrm{Id}+\cL_{\cE_\bS}-\nabla\fE,\
\Phi_\xi(\eta)=-(\xi\star\eta),\ R_0=\fE\star{}=-\Phi_\fE,\ q=0,\ w=2-D.
\]

\subsection{$\bS$-Equivariant Frobenius manifolds with finite group
action}\label{finite group} Let us consider an $\bS$-equivariant
Frobenius manifold $(M,\star,g,e,\fE)$ of weight $D$ over $R$. Let
$W$ be a finite group which acts by automorphisms on $M$ (hence on
$\Theta_M$) in a compatible way with the $\bS$-equivariant
Frobenius structure. In particular, the action of $W$ on
$\Theta_M$ commutes with $\cL_{\cE_\bS}$.

Let $M^W$ be the fixed set of $W$ on $M$. Then $W$ acts by $\cO_{M^W}$-linear
isomorphisms on $\Theta_M\big|_{M^W}$. Moreover, the fixed set $M^W$ is a
smooth subscheme of $M$ over $R$ and the fixed bundle $(\Theta_M\big|_{M^W})^W$
is equal to $\Theta_{M^W}$.

Let us moreover assume that $W$ is equipped with a \emph{non trivial} character
$\sgn:W\to\nobreak\{\pm1\}$. We denote by $a :\Theta_M\big|_{M^W}\to
\Theta_M\big|_{M^W}$ the antisymmetrization morphism and by $E$ its image. Then
$E$ is a locally free $\cO_{M^W}$-submodule of $\Theta_M\big|_{M^W}$ and we
have a decomposition $\Theta_M\big|_{M^W}=E\oplus\ker a $. This decomposition
is $g$-orthogonal, as $g(a\xi,a\eta)=g(\xi,\eta)$ for any $\xi,\eta$ and $g$
restricted to $E$ is nondegenerate.

As the inclusion $M^W\hookrightarrow M$ commutes with $\cL_{\cE_\bS}$, one can
restrict to $M^W$ the $\bS$-equivariant pre-Saito structure associated to
$(M,\star,g,\fE)$ to get such an object with corresponding bundle
$\Theta_{M|M^W}$. One can moreover induce this structure on the
$\cO_{M^W}$-module $E$, as the following operators leave $E$ invariant:
\begin{itemize}
\item
the connection $\nabla$ (i.e., $\nabla_\xi\eta$ is a section of $E$ whenever
$\xi$ is a section of $\Theta_{M^W}$ and $\eta$ a section of $E$),
\item
the Higgs field $\Phi$, (i.e., $\xi\star\eta=-\Phi_\xi\eta$ is a section of $E$
whenever $\xi$ is a section of $\Theta_{M^W}$ and $\eta$ a section of $E$),
\item
the operator $R_0=-\Phi_\fE=\fE\star{}$,
\item
the operator $\wt\nabla_{\cE_\bS}=\mathrm{Id}+\cL_{\cE_\bS}-\nabla\fE$ (i.e.,
$\nabla_\eta\fE$ is a section of $E$ whenever $\eta$ is a section of $E$).
\end{itemize}

The following is then clear:
\begin{Lemma}
The tuple $(M^W,E,\wt\nabla,\Phi,R_0,g)$ defines an
$\bS$-equivariant pre-Saito structure of weight $w=2-D$ on $M^W$.
\end{Lemma}

\begin{Prop}\label{prop:inducedFrobenius}
Let us assume that there exists $\omega\in E\subset\Theta_{M|M^W}$
which is $\nabla$-horizontal and an eigenvector of
$\wt\nabla_{\cE_\bS}$ (acting on $E$ or on $\Theta_{M|M^W}$) and
such that the morphism
\begin{align*}
\Theta_{M^W}&\longrightarrow E\\
\xi&\longmapsto\xi\star\omega
\end{align*}
is onto. Then, any smooth formal subscheme $N\subset M^W$ over $R$
defined by an ideal invariant under $\cL_{\cE_\bS}$ and such that
the induced morphism $\Theta_N\to E_{|N}$ is an isomorphism comes
equipped with a natural structure of an $\bS$-equivariant
Frobenius manifold of weight $D$.
\end{Prop}

\begin{proof}
We restrict the $\bS$-equivariant pre-Saito structure
$(M^W,E,\wt\nabla,\Phi,R_0,g)$ to $N$ to get an object
$(N,E_{|N},\wt\nabla,\Phi,R_0,g)$ of the same kind. Then, as $\omega_{|N}$ is a
$\nabla$-horizontal section of $E_{|N}$ and as the morphism $\Theta_N\to
E_{|N}$ given by $\xi\mapsto\xi\star\omega_{|N}=\varphi_\omega(\xi)$ is an
isomorphism, $\omega$ is primitive. Moreover, $\omega$ is homogeneous in $E$
hence $\omega_{|N}$ is so in $E_{|N}$. One can then apply the correspondence
of~\S\ref{equiv Frobenius}.
\end{proof}

\subsubsection*{Some properties of the $\bS$-equivariant Frobenius manifold structure on $N$}
Abusing notation, we denote by ${}\star\omega^{-1}$ the inverse map of the
induced isomorphism ${}\star\omega:\Theta_N\to E_{|N}$

We denote by ${}^\omega\!g,\onabla$ the metric and connection on $\Theta_N$
coming from that on $E_{|N}$, and by $\circ$ the product on $\Theta_N$ induced
by the Higgs field on $E_{|N}$.
\begin{enumerate}
\item\label{item:star}
Let $\xi,\eta$ be sections of $\Theta_N$. The product $\xi\star\eta$ in
$\Theta_{M|N}$ may not be a section of $\Theta_N$ (it is only a section of
$\Theta_{M^W|N}$). We have $[\xi\star\eta-\xi\circ\eta]\star\omega=0$. In fact,
the composition
\[
\Theta_{M^W|N}\stackrel{\star\omega}\longrightarrow
E_{|N}\stackrel{\star\omega^{-1}}\longrightarrow \Theta_N
\]
induces a projection $\Theta_{M^W|N}\to \Theta_N$, and $\xi\circ\eta$ is
nothing but the projection of $\xi\star\eta$ on $\Theta_N$, a formula that we
can write
\[
\xi\circ\eta=(\xi\star\eta\star\omega)\star\omega^{-1}.
\]

\item\label{item:unit}
Let us assume that we can find $N$ such that the unit field $e$ is
\emph{tangent} to $N$. This condition does not lead to a contradiction, as
$e\star\omega=\omega\neq0$. Then~$e_{|N}$ is the unit field for the
$\bS$-equivariant Frobenius manifold structure on $N$. Indeed, clearly,
$e_{|N}\circ\eta=\eta$ for any section $\eta$ of $\Theta_N$. On the other hand,
we have to check that $e$ is $\onabla$-horizontal:
\[
\onabla e_{|N}:=\nabla
(e_{|N}\star\omega)\star\omega^{-1}=\nabla(\omega)\star\omega^{-1}=0,
\quad\text{as }\nabla(\omega)=0.
\]

\item\label{item:Euler}
Let us assume that $N$ is chosen so that the Euler vector field $\fE$ is
tangent to $N$. Then $\fE_{|N}$ is the Euler vector field for the Frobenius
manifold structure on $N$, as $R_0=\fE\star$ leaves $E$ invariant.

\item\label{item:metric}
We have ${}^\omega\!g(\xi,\eta)=g(\xi\star\omega,\eta\star\omega)$ for any
$\xi,\eta\in\Theta_N$.
\end{enumerate}

\begin{Rmk}
Given an $R$-basis $\boldsymbol{e}^o$ of $E/(t_0,\dots,t_{m-1})E$,
there exists a unique system of flat coordinates $(t_i)$ on $N$
for which $\partial_{t_i}\star\omega\equiv e_i^o
\bmod(t_0,\dots,t_{m-1})E$. Given any other formal smooth
subscheme $N'$ over $R$ satisfying the properties in Proposition
\ref{prop:inducedFrobenius}, with corresponding coordinates
$(t'_i)$, we do not know whether the natural isomorphism
$\cO_N\to\cO_{N'}$, $t_i\mapsto t'_i$, is compatible with the
$\bS$-equivariant Frobenius structures. In other words, there is a
priori no uniqueness in the construction resulting from
Proposition \ref{prop:inducedFrobenius}. However, when this
construction is applied to the setting of \S\ref{Setting},
Conjecture \ref{conj} also gives uniqueness.
\end{Rmk}

\section{The abelian/nonabelian correspondence for Frobenius structures}\label{Frobenius correspondence}

A precise relation between the genus zero Gromov-Witten theory
(with descendants) of a quotient by a nonabelian group and a twist
of the theory for the quotient by a maximal torus in the group was
conjectured in \cite{BCK2}. Here we formulate a version of this
correspondence for the associated Frobenius structures.

\subsection{Setting}\label{Setting} Let $X$ be a
smooth projective variety over
${\CC}$ with the (linearized) action of a complex reductive group $\G$, and let
$\T \subset \G$ be a maximal torus. In this setting, there are
two geometric invariant theory (GIT) quotients, $X//\T$ and $X//\G$.
We assume (for both actions) that all semistable points are stable and that
all isotropy groups of stable points are trivial, so that
$X//\T$ and $X//\G$ are smooth projective varieties. Further, we assume
that the $\G$-unstable locus $X\setminus X^s(\G)$ has codimension at least 2 in $X$.
(Note that this last condition is automatic when X is a projective space.)

There is a diagram

\[\begin{CD} X//\T = X^s(\T)/\T @<<{j}< X^s(\G)/\T \\
@. @VV{\pi}V \\
@. X//\G=X^s(\G)/\G
\end{CD}\]
with $j$ an open immersion and $\pi$ a $\G/\T$-fibration.

The above diagram leads to a comparison of the cohomology of the nonabelian
quotient $X//\G$ to that of the abelian quotient $X//\T$ (\cite{ES},\cite{Martin},\cite{Kir}).
We describe an equivariant version of it.
Let another (possibly trivial) complex torus
$\bS$ act on $X$. Assume
that the action commutes with the action of $\G$ and preserves
$X^s(\T)$ and $X^s(\G)$. There is an induced action of $\bS$ on
the smooth projective varieties $X//\T$ and $X//\G$. The morphisms
in the diagram are $\bS$-equivariant. To the pair $(\G,\T)$
we associate the usual Lie-theoretic data:
\begin{itemize}
\item the Weyl group
$\bW=N(\T)/\T$ ($N(\T)$ is the normalizer of $\T$ in
$\G$).
\item the root system $\Phi$ with decomposition $\Phi=\Phi_+\cup\Phi_{-}$
into positive and negative roots.
\item for each root $\alpha\in \Phi$ the $1$-dimensional $\T$-representation
$\CC_{\alpha}$ with weight $\alpha$.
\end{itemize}
The Weyl group acts on $X//\T$, hence also on the equivariant cohomology ring
$H^*_{\bS}(X//\T,\CC)$.
The representations $\CC_\alpha$ define $\bS$-equivariant line bundles
$$L_\alpha := X^s(\T)\times _\T \CC
_{\alpha}$$ on $X//\T$ with equivariant first Chern classes
$c_1^{\bS}(L_\alpha)\in H^*_{\bS}(X//\T,\CC)$.
The $\bS$-action on $L_\alpha$ is induced
by the $\bS$-action on $X^s(\T)$ (and the trivial $\bS$- action on
$\CC _\alpha$). Note that $L_{-\alpha}\cong L_{\alpha}^{\vee}$ for
any pair $(\alpha,-\alpha)$ of opposite roots. The equivariant cohomology class
$$ \omega := \sqrt{\frac{1}{|\bW|}\prod _{\alpha\in\Phi}
c_1^{\bS}(L_\alpha )}=\sqrt{\frac{(-1)^{|\Phi_+|}}{|\bW|}}
\prod_{\alpha\in\Phi_+}c_1^{\bS}(L_\alpha ) $$
will play an important role in this paper. It is the fundamental
$\bW$-anti-invariant class in the equivariant cohomology of $X//\T$; any
other $\bW$-anti-invariant class $\phi$ can be written
(non-uniquely) as $\gamma\cup\omega$,
with $\gamma\in H^*_\bS(X//\T,\CC)^\bW$.
(The reason for considering $\omega$
rather than the customary
$\Delta=\prod_{\alpha\in\Phi_+}c_1^{\bS}(L_\alpha )$
is one of convenience: we simply want to avoid having to insert the factor
${(-1)^{|\Phi_+|}}/{|\bW|}$ in all formulae comparing Gromov-Witten
invariants of $X//\G$ and $X//\T$.)

The following facts are known:

$(3.1.1)$ $\pi ^*$ induces an isomorphism $\
H^*_{\bS}(X//\G)\cong H^*_{\bS}(X^s(\G)/\T)^\bW$

$(3.1.2)$ There is an exact sequence
\[ \begin{array}{ccccccccc}
0 &\longrightarrow & \mathrm{ker}(\cup\omega ) &\longrightarrow & H^*_{\bS}(X//\T)^\bW &
\longrightarrow & H^*_{\bS}(X//\G) &
\longrightarrow & 0 \\
& & &\subset& & (\pi ^*)^{-1}\circ j^* & & &
\end{array} \] where $\mathrm{ker}(\cup\omega )$ is $\{ \gamma \in
H^*_{\bS}(X//\T)^\bW \ | \ \gamma \cup \omega =0 \}$.

$(3.1.3)$ The equivariant push-forwards satisfy the equality
\[ \int _{X//\T} \omega ^2 \tilde{\sigma} = \int
_{X//\G} \sigma \] for all $\sigma \in H^*_{\bS}(X//\G)$,
$\tilde{\sigma}\in H^*_{\bS}(X//\T)$ with
$j^*\tilde{\sigma}=\pi ^*(\sigma )$. (Such $\tilde{\sigma}$ are
called {\em lifts} of $\sigma$.)

$(3.1.4)$ There is an identification of the $\bS$-equivariant
relative tangent bundle $T_{\pi}$ of $\pi : X^s(\G)/\T \ra
X^s(\G)/\G$ with $\oplus_{\alpha \in\Phi} L_\alpha
|_{X^s(\G)/\T}$.

\medskip

In the nonequivariant case (that is, for $\bS$=1),
$(3.1.1)$ is a classical result, $(3.1.2)$ is proved in \cite{ES} for $X=\PP^N$ and
in \cite{Kir} in general,
%$\footnote{By
%Theorem A in \cite{Martin} and lines 3-4 on page 171 in \cite{ES},
%in fact, $\mathrm{ann}(\omega \cup ) = \mathrm{ann}(\omega ^2\cup
%)$},
$(3.1.3)$ is proved in \cite{Martin} and $(3.1.4)$ can be seen by a direct
computation. The extensions to the equivariant context are straightforward
and left to the reader.

\subsection{The $\bW$-induced Frobenius manifold}\label{induced}
Applying the results in \S \ref{finite group} to the Weyl group
action on the $\bS$-equivariant Frobenius manifold given by the
equivariant Gromov-Witten theory of $X//\T$, a new
$\bS$-equivariant Frobenius manifold (of dimension over the base
ring equal to the rank of $H^*_{\bS}(X//\G,\CC)$) is obtained. In
this subsection we spell out for concreteness the details of the
construction and the main properties of the new Frobenius
structure in this special case.

As mentioned in the introduction, a specialization of Novikov variables will be
needed before comparing the new Frobenius structure with the one given by the
equivariant Gromov-Witten theory of $X//\G$ and we start with this specialization.

Recall from (3.1.3) the notion of lift of cohomology classes from
$X//\G$ to $X//\T$. By $(3.1.2)$, one can always choose
$\bW$-invariant lifts. These are not generally unique; however,
the assumption that the $\G$-unstable locus in $X$ has codimension
$\geq 2$ implies that for divisor classes the $\bW$-invariant
lifts are unique. This allows us to lift curve classes as well
(cf. \cite{BCK2}): the inclusion
$$\Pic(X//\G)\cong \Pic(X//\T)^{\bW}
\subset \Pic(X//\T)$$
induces by duality a projection
$$\varrho :NE_1(X//\T)\longrightarrow NE_1(X//\G).$$
We say that $\tilde{\beta}$ lifts $\beta\in NE_1(X//\G)$
(and write $\tilde{\beta}\mapsto\beta$) if $\varrho(\tilde{\beta})=\beta$.
Note that any effective $\beta$ has finitely many lifts.
Define a projection on Novikov rings
\begin{equation}\label{specialization} p:N(X//\T)\ra N(X//\G),
\;\;\;\;\; p\left(\sum_{\tilde{\beta}}c_{\tilde{\beta}}Q^{\tilde{\beta}}\right)=
\sum_\beta(-1)^{\epsilon(\beta)}
\left(\sum_{\tilde{\beta}\mapsto\beta}c_{\tilde{\beta}}\right)Q^\beta ,
\end{equation}
where $$\epsilon :NE_1(X//\G)\longrightarrow \ZZ_2$$
is defined by
$$\epsilon(\beta)=
\left(\int_{\tilde{\beta}}\sum_{\alpha\in\Phi_+}c_1^{\bS}(L_\alpha)\right)\;({\mathrm {mod}} 2)$$
with $\tilde{\beta}$ any lift of $\beta$. This makes sense, since the right-hand side
does not depend on the choice of lift. Indeed, if $\alpha'$ is any simple root and $v_{\alpha'}\in\bW$
is the corresponding reflection, then by standard properties of root systems
$$ v_{\alpha'}(\sum_{\alpha\in\Phi_+}c_1^{\bS}(L_\alpha))=\sum_{\alpha\in\Phi_+}c_1^{\bS}(L_\alpha)-
2c_1^\bS(L_{\alpha'}),$$
so $\sum_{\alpha\in\Phi_+}c_1^{\bS}(L_\alpha)$
is $\bW$-invariant as a cohomology class with $\ZZ_2$-coefficients.

The sign in (\ref{specialization}), which may seem rather mysterious, has its
origin in the twisting bundle appearing in the abelian/nonabelian correspondence,
as formulated in Conjecture 4.2 of \cite{BCK2}.

Let $Z$ be the formal Frobenius manifold defined by the $\bS$-equivariant
Gromov-Witten theory of $X//\T$, with potential function $F^{X//\T,\bS}$. Let $M$ be the formal scheme over
$N(X//\G)\otimes\CC[\lambda]$ obtained by base change from $Z$ by
the morphism of Novikov rings (\ref{specialization}). Let $\theta:M\lra Z$ be the
base change map. We obtain a formal Frobenius structure over $N(X//\G)\otimes\CC[\lambda]$
on $(M, \Theta_M)$
by pulling-back via $\theta$ the Frobenius structure on $Z$. Note that only
the potential (and therefore the quantum product) changes under the pull-back,
since the coefficients of the metric, the horizontal sections and the Euler vector field
do not depend on the Novikov variables. Explicitly, the potential of the Frobenius
structure on $M$ is
\begin{equation}\label{potential}
F:=\theta^*(F^{X//\T,\bS})=\sum_{\beta\in NE_1(X//\G)}(-1)^{\epsilon(\beta)}Q^\beta\sum_{n\geq 0}\frac{1}{n!}
\left(\sum_{\tilde{\beta}\mapsto\beta}
\langle\underbrace{\gamma,\dots, \gamma}_{n} \rangle_{0,n,\tilde{\beta}}^{X//\T,\bS}\right)
\end{equation}

Choose a homogeneous basis $\{ \sigma_0=1,\sigma_1\dots ,\sigma_r,\sigma_{r+1}\dots,
\sigma_{m-1} \} $ of $H^*_{\bS}(X//\G)$ over
$\CC[\lambda]:=\CC[\lambda_1,\dots ,\lambda_{\ell}]=H^*(B\bS)$, such that
$\{\sigma_1,\dots ,\sigma_r\}$ forms a basis of $H^2_{\bS}(X//\G)$
and {\underline {fix}}
homogeneous lifts $\gamma _i \in H^*_{\bS}(X//\T)^\bW$ of $\sigma _i$.
The fixed lifts give rise to a $\CC$-linear embedding
\begin{equation}\label{embedding} H^*_{\bS}(X//\G,\CC )\subset
H^*_{\bS}(X//\T,\CC )\end{equation} (which may not in general be a
homomorphism of equivariant cohomology rings).

The image of the embedding (\ref{embedding}) determines
a formal submanifold $N$ of $M$ over $N(X//\T)\otimes\CC [\lambda]$.

Let $$V:=H^*_{\bS}(X//\T,\CC)^a$$ be the subspace of $\bW$-anti-invariant
classes. The composition of (\ref{embedding}) with the map
$$\cup\omega: H^*_{\bS}(X//\T,\CC)^\bW\ra H^*_{\bS}(X//\T,\CC)^a$$ is an isomorphism from
$H^*_{\bS}(X//\G,\CC)$ to $V$. Let $\cV=V\otimes \cO_N$ be the
subsheaf of $\Theta_M|_N$ induced by $V$. Let $\star$ be the
quantum product on $\Theta_M$ (that is, the pull-back by $\theta$
of the quantum product on $H^*_\bS(X//\T,\CC)$). Consider the map
$$\star\omega:\left(\Theta_M|_N\right)^{\bW}
\longrightarrow \cV,\;\;\;\; \xi\mapsto(\hat{\xi}\star\omega)|_N,$$
with $\hat{\xi}\in \Theta_M^{\bW}$ any extension of $\xi$ to $M$.
(It is well defined, since the quantum product of two vector
fields at a point depends only on their values at the point.) This
map reduces to $\cup\omega$ modulo the ideal generated by
$\{Q^{\beta}|\beta\neq 0\}$. By Nakayama's Lemma, $\star\omega$
induces an isomorphism $\Theta_N\ra\cV$. Let $\phi:\cV\ra
\Theta_N$ be the inverse isomorphism. Abusing notation, when
$\eta\in\cV$ we write $\eta\star\omega^{-1}$ for $\phi(\eta)$.
Hence we have for $\xi \in \Theta_N$
$$(\xi \star\omega)\star\omega^{-1}=\xi .$$

We now induce a structure of formal Frobenius manifold on $N$
(over $N(X//\G)\otimes\CC [\lambda]$) using the maps $\star\omega$
and $\star\omega^{-1}$. Explicitly:

\medskip
\noindent (3.2.4) The metric $^{\omega}g$ on $\Theta_N$ is given by the composition
$$
\begin{CD} \Theta_N\otimes \Theta_N\hookrightarrow
\Theta_M|_N\otimes \Theta_M|_N@>\star\omega\otimes
\star\omega>>\cV\otimes\cV@>g|_{\cV}>>\cO_N,
\end{CD}$$
that is,
$$^{\omega}g(\xi,\eta)=g|_{\cV}(\xi \star\omega, \eta \star\omega).$$
Note that $g|_{\cV}$ is nondegenerate on $\cV$ by Martin's formula
(3.1.3).

\medskip
\noindent (3.2.5) The Levi-Civita connection $^{\omega}\nabla$ of $^\omega g$ satisfies
$$^{\omega}\nabla_{\xi}\eta=\left(\nabla_{\hat{\xi}}(\hat{\eta} \star\omega)\right)|_N\star\omega^{-1}.$$
%(To see this, it suffices to check that the right-hand side defines a connection which is
%torsion-free and compatible with the metric $^\omega g$. These are routine calculations and are omitted.)
\medskip

\noindent (3.2.6) The product of $\xi,\eta\in \Theta_N$ is defined by
$$\xi \circ \eta = (\xi \star \eta \star
\omega ) \star\omega ^{-1}.
$$
In other words, $\xi\circ\eta$ is the projection of $\xi
\star\eta$ along $\ker (\star\omega)$.

\medskip

\noindent(3.2.7) The unit is the vector field $1$ restricted to $N$.

\medskip

The symmetry of $^\omega\nabla( ^\omega g(\cdot \circ\cdot ,\cdot ))$
and the  corresponding potential function
are discussed in $\S$ 3.5 below.

%The reader can easily check that all other requirements in Definition \ref{Frob} are
%satisfied.

\subsection{Flat coordinates}
On $N$ there are coordinates $\tilde{t}_0,...,\tilde{t}_{m-1}$
determined by the basis $\{\gamma _0 =1, \gamma _1,...,\gamma _r,...,\gamma
_{m-1}\}$ of lifts introduced above. These are just restrictions to $N$
of coordinates on $M$ which are flat for the connection $\nabla$.
Let
\begin{equation} \label{horizontal}
\xi_i(\tilde{t}):=(\gamma_i\cup\omega)\star\omega^{-1},\;\;\;\;i=0,\dots,{m-1}.\end{equation}
Equivalently, $\xi_i$ is defined by the equality
\begin{equation}\label{flat}
\xi_i\star\omega=\gamma_i\cup\omega .
\end{equation}
The $\xi_i$'s form a basis of $\Theta_N$
consisting of $^\omega\nabla$-horizontal vector fields.
Denote by $s:=(s_0,s_1,...,s_r,...,s_{m-1})$ the corresponding
$^\omega\nabla$-flat coordinates
on $N$ (so that $\partial_{s_i}=\xi_i$). Note that
\begin{equation}\label{first order}
s \equiv \tilde{t}, \text{ modulo the ideal generated by }
 \{Q^{{\beta}}|{\beta}\neq 0\}.
\end{equation}

\subsection{Euler vector field}\label{degree}
Since
\begin{eqnarray*} j^*c_1^{\bS}(T(X//\T)) &=& c_1^{\bS} (T(X^s(\G)/\T)) \\
&=& \pi ^*(c_1^{\bS}(T(X//\G))) + \sum _{\alpha\in\Phi} c_1^{\bS}(L_\alpha ) \\
&=& \pi ^* (c_1^{\bS}(T(X//\G)))\end{eqnarray*}
and $ \Pic _{\bS} (X//\T) \cong
\Pic _{\bS}(X^s(\G)/\T)$ via $j^*$, we conclude that $c_1^{\bS}(T(X//\T))$
is $\bW$-invariant. Viewing $c_1^{\bS}(T(X//\T))$ as a
vector field on $M$, its restriction to $N$ is therefore a section of $(\Theta_M|_N)^{\bW}$.
Moreover,
this restriction is in fact tangent to $N$, since (by uniqueness of lifts
of divisors)
$N$ contains the germ of linear subspace $H^2_{\bS}(X//\T)^\bW$.

Define the Euler vector field by
$$\fE_\bS = \fE + \sum _{i=1}^\ell \lambda _i
\partial _{\lambda _i}$$
with
$$\fE = \sum _{i=0}^{m-1} (1-\frac{\mathrm{ deg} \sigma
_i}{2})\tilde{t}_i
\partial _{\tilde{t}_i} + c_1^{\bS}(T(X//\T))|_N.$$
Note that $\fE_\bS$ is simply the restriction to $N$ of the
corresponding Euler vector field for $X//\T$ (see \S 2.2): Extend
$\{\gamma _0,...,\gamma _{m-1}\}$ to a basis of $H^*_{\bS}(X//\T)$
for the Euler vector field for $X//\T$.

Applying $\cL_{\fE_\bS}$ to the equality $(\xi_i\star\omega ) =
\partial _{\tilde{t}_i}\cup \omega$, we see that
$$\cL _{\fE_\bS}
\xi_i = (1-\frac{\text{ deg } \sigma
_i}{2})\xi_i.$$
Easy calculations show then that
$$\cL _{\fE_\bS}(^\omega g)=(2-\dim(X//\G))\; ^\omega g,\;\;\;\;
\cL _{\fE_\bS}(\circ )=\circ$$ (and obviously $\cL
_{\fE_\bS}(1)=-1$), hence $\fE_{\bS}$ is indeed an Euler vector
field. Also,
$$\cL _{\fE_\bS} (s_i) = \deg
(\tilde{t}_i) s_i,$$
that is, $\deg s_i = \deg \tilde{t}_i$. In particular, $\deg
s_1=...=\deg s_r=0$.

\subsection{Potential}
Recall that we identify the cohomology classes on $X//\T$ with
$\cO_M$-linear vector fields on $M$. Denote by $
\partial _{\tilde{t}_i\cup \omega}$ the vector field
corresponding to $\gamma_i\cup\omega$. The components of the
tensor \lq\lq $\circ$" in the basis of $^\omega\nabla$-horizontal
fields are
\begin{equation} \label{3.5.1} ^\omega g (\xi_i\circ
\xi_j,\xi_k) = g|_{\cV} (\xi_i\star \xi_j \star \omega ,\xi_k
\star\omega )=
g(\hat{\xi_i}\star(\gamma_j\cup\omega),\gamma_k\cup\omega)|_N
\end{equation}
where $\hat{\xi_i}$ is any extension of $\xi$ to a $\bW$-invariant vector field on $M$.
Since
\begin{multline}\notag
g(\hat{\xi_i}\star(\gamma_j\cup\omega),\gamma_k\cup\omega)|_N=
g(\gamma_j\cup\omega,\hat{\xi_i}\star(\gamma_k\cup\omega))|_N=\\
g|_{\cV}(\gamma_j\star\omega,\xi_i\star\xi_k\star\omega)= ^\omega
g(\xi_j,\xi_i\circ\xi_k),
\end{multline}
we see that the Frobenius algebra property
\begin{equation}\label{frobeniusalgebra} ^\omega g (\xi_i\circ
\xi_j,\xi_k) = {^\omega} g(\xi_j,\xi_i\circ\xi_k)
\end{equation}
holds.
Recall the potential $F$ (see (\ref{potential})) of the formal
Frobenius manifold $M$. We get from (\ref{3.5.1})
\begin{multline}\label{new product}
^\omega g (\xi_i\circ
\xi_j,\xi_k) =g(\hat{\xi_i}\star(\gamma_j\cup\omega),\gamma_k\cup\omega)|_N=\\
\left.\left({{\hat{\xi}}_i}
(\partial _{\tilde{t}_j\cup \omega}
\partial _{\tilde{t}_k\cup \omega} F )\right)\right|_N=
\xi_i\left(\left. (\partial _{\tilde{t}_j\cup \omega}
\partial _{\tilde{t}_k\cup \omega} F )\right|_N\right)  .\end{multline}
Note that
$$
%\xi_l\left(
%\hat{\xi}_i(\partial _{\tilde{t}_j\cup \omega } \partial _{\tilde{t}_k\cup
%\omega } F)) |_N\right)=
\xi_l(\xi_i((\partial _{\tilde{t}_j\cup \omega } \partial _{\tilde{t}_k\cup
\omega } F) |_N))=
\xi_i(\xi_l((\partial _{\tilde{t}_j\cup \omega } \partial _{\tilde{t}_k\cup
\omega } F) |_N))
%=\xi_i\left(
%\hat{\xi}_l(\partial _{\tilde{t}_j\cup \omega } \partial _{\tilde{t}_k\cup
%\omega } F)) |_N\right)
,$$
since $[\xi_l,\xi_i]=0$.
Hence
$$\xi_l(^\omega g (\xi_i\circ
\xi_j,\xi_k))=\xi_i(^\omega g (\xi_l\circ
\xi_j,\xi_k)).$$
Combined with (\ref{frobeniusalgebra}), this implies that the tensor
$\xi_l(^\omega g (\xi_i\circ
\xi_j,\xi_k))$
is symmetric in the indices $l,i,j,k$, hence
there is a (formal) function $F'$ on $N$ such that $$ \partial
_{s_i}\partial _{s_j}
\partial _{s_k} F' = {^\omega g}(\xi_i\circ \xi_j,\xi_k). $$
Finally, a direct computation shows that $L_{\fE
_\bS} F' = (3 - \dim X//\G ) F'$ up to quadratic terms.

This finishes the construction of the induced structure of formal
$\bS$-equivariant Frobenius manifold over
$N(X//\G)\otimes\CC[\lambda]$ on $N$.

\subsection{More on the flat coordinates} For later use we record here
some details about the
change of coordinates $s_i(\tilde{t})$ on $N$.
From the defining equation (\ref{horizontal}) for the horizontal vector fields
$\xi_i$ it follows that the jacobian matrix $A:=(\partial s_i/\partial\tilde{t}_j)_{i,j}$ is
given explicitly by
\begin{equation*} \frac{\partial s_i}{\partial\tilde{t}_j}=\left.(\partial_{\tilde{t}_j}
\partial_{\tilde{t}_i\cup\omega}^{\vee}\partial_{\tilde{t}_0\cup\omega}F)\right|_N
\end{equation*}
where
$$\partial_{\tilde{t}_i\cup\omega}^{\vee}:=\sum_k g^{ik}\partial_{\tilde{t}_k\cup\omega}
$$
with $(g^{ik})\in GL_m(\CC[\lambda])$ the inverse matrix of the metric $g|_\cV$.
Using the divisor axiom for Gromov-Witten invariants of $X//\T$ in the formula
(\ref{potential}) for the potential $F$, we see that the entries of the jacobian matrix
have the form
\begin{eqnarray}\label{jacobian} \notag\frac{\partial s_i}{\partial\tilde{t}_0}&=&\delta_{i0}\\
\frac{\partial s_i}{\partial\tilde{t}_j}&=&\delta_{ij}+\sum_{\beta\neq 0}Q^{\beta}
e^{\beta\cdot \tilde{t}_{\mathrm{small}}}c_{\beta,ij}(\tilde{t}_{r+1},\dots,\tilde{t}_{m-1}),\;\;\;\;\; j\neq 0,
\end{eqnarray}
where $c_{\beta,ij}\in\CC[\lambda][[\tilde{t}_{r+1},\dots,\tilde{t}_{m-1}]]$ and
$$\beta\cdot\tilde{t}_{\mathrm{small}}=\sum_{i=1}^r\tilde{t}_i\int_{\beta}\sigma_i.$$
By integrating (\ref{jacobian}) (with the initial condition $s(0)=0$),
we obtain a refined version of (\ref{first order})
\begin{equation}\label{mirrortransf}
%s_0&=&\tilde{t}_0\\
s_i=\tilde{t}_i+\sum_{\beta\neq 0}Q^{\beta}
e^{\beta\cdot \tilde{t}_{\mathrm{small}}}b_{\beta,i}(
\tilde{t}_{r+1},\dots,\tilde{t}_{m-1}),
\end{equation}
with $ b_{\beta,i}\in\CC[\lambda][[\tilde{t}_{r+1},\dots,\tilde{t}_{m-1}]]$.

By considering the inverse jacobian matrix (which gives the map
$\star\omega^{-1}$), it follows that the inverse coordinate change
$\tilde{t}(s)$ is also of the type (\ref{mirrortransf}), hence the
potential function $F'$ in flat coordinates $s_i$ has the special
form (\ref{qc-type}) (up to $\leq$ quadratic terms in the $s_i$'s)
\begin{equation}\label{qc-type'} F'=F'_{cl}+\sum_{\beta\neq 0}Q^{\beta}
e^{\beta\cdot s_{\mathrm{small}}}F'_{\beta}(s_{r+1},\dots,s_{m-1})
\end{equation}
where $ \beta\cdot s_{\mathrm{small}}=\sum_{i=1}^r s_i(\int_{\beta}\sigma_i)$.

Finally, we record what happens with the \lq\lq small" parameter
spaces under the change of coordinates.
\begin{Lemma}\label{small space} $(i)$ If $X//\G$ is Fano of index $\geq 2$,
then the subspaces of $N$ given by the equations $\{s_0=s_{r+1}=\dots =s_{m-1}=0\}$ and
$\{\tilde{t}_0=\tilde{t}_{r+1}=\dots =\tilde{t}_{m-1}=0\}$ coincide. Moreover,
on this subspace we have $s_i=\tilde{t}_i$ for $i=1,\dots, r$.

$(ii)$ If $c_1(T(X//\G)$ is nef, then the subspaces $\{ s_{r+1}=\dots =s_{m-1}=0\}$
and $\{\tilde{t}_{r+1}=\dots =\tilde{t}_{m-1}=0 \}$ coincide.

\end{Lemma}

\begin{proof} $(i)$ Let $1\leq i\leq r$. After restriction to
$\tilde{t}_0=\tilde{t}_{r+1}=\dots=\tilde{t}_{m-1}=0$ we obtain
$$\xi_i=\gamma_i+\sum_j\left(\sum_{\beta}c_{\beta,ij}
Q^\beta e^{\beta\cdot
\tilde{t}_{\mathrm{small}}}\right)\gamma_j.$$ Since
$\deg{\xi_i}=1$ and $\deg{e^{\beta\cdot
\tilde{t}_{\mathrm{small}}}}= \int_{\beta}c_1(T(X//\G)\geq 2$, we
must have $\xi=\gamma_i$ and the statement follows. The proof of
$(ii)$ is similar. \end{proof}

\subsection{Main Conjecture}
Let $P$ be the formal $\bS$-equivariant Frobenius manifold over
$N(X//\G)\otimes\CC[\lambda]$ defined by the genus zero
$\bS$-equivariant Gromov-Witten theory of $X//\G$, with flat
coordinates $t_0,t_1,\dots ,t_r,\dots,t_{m-1}$ corresponding to
the $\CC [\lambda]$-basis
$\{\sigma_0=1,\sigma_1,\dots,\sigma_r,\dots,\sigma_{m-1}\}$ of
$H^*_{\bS}(X//\G)$ and potential function $F^{X//\G,\bS}$. We are
now ready to formulate the abelian/nonabelian correspondence:
\begin{Conj} \label{conj} Let $\varphi :P\longrightarrow N$ be the isomorphism of
formal schemes over $N(X//\G)\otimes\CC[\lambda]$ defined by
$\varphi^*(s_i)=t_i$. Then $\varphi$ induces an isomorphism of
formal $\bS$-equivariant Frobenius structures such that
$\varphi^*(\xi_i)=\sigma_i$ and $\varphi^*F'=F^{X//\G,\bS}$ up to
quadratic terms.
\end{Conj}

Note that $\varphi^*(\xi_i)=\sigma_i$ follows easily from (\ref{horizontal}). The main
point of the conjecture is the identification of potentials. We also remark
that the conjecture implies in particular that the new $\bW$-induced
Frobenius structure constructed in this section does not depend on the
choice of the $\bW$-invariant lift of $H^*_{\bS}(X//\G,\CC)$.

\section{Proof of Conjecture \ref{conj} for flag manifolds}\label{the proof}
\subsection{Preliminaries}
Let $0<k_1<\dots <k_r<n=k_{r+1}$ be integers. Consider the
vector space
$$\Omega :=\bigoplus_{i=1}^r{\mathrm{Mat}}_{k_i\times k_{i+1}}(\CC)$$
where ${\mathrm{Mat}}_{k_i\times k_{i+1}}(\CC)$ is the space of
matrices of size $k_i\times k_{i+1}$ with complex entries.
Let $\G:=\prod_{i=1}^r GL_{k_i}(\CC)$, with maximal torus $\T$ equal to the
product of the subgroups of diagonal matrices. $\G$ acts on $\Omega$ by
$$(g_1,\dots ,g_r)(A_1,\dots A_r)=(g_1A_1g_2^{-1},g_2A_2g_3^{-1},\dots,
g_{r-1}A_{r-1}g_r^{-1},g_rA_r).$$ This action descends to an
action on $X:=\PP(\Omega)$, with a canonical linearisation on
$\cO(1)$, and the GIT quotient $X//\G$ is the partial flag
manifold $Fl(k_1,\dots,k_r,n)$ parameterizing flags of subspaces
$\{\CC^{k_1}\subset\dots\subset \CC^{k_r}\subset\CC^n\}$.

The corresponding abelian quotient $X//\T$ is a toric variety which can be
realized as a tower of fibered products of projective bundles.

Let $\bS\cong(\CC^*)^n$ be the subgroup of diagonal matrices in $GL_n(\CC)$,
acting on $\Omega$ by right-multiplication of $A_r$. There are induced $\bS$-actions
on $X//\G$ (which is just the usual action of the maximal torus in $GL_n$ on
the flag manifold) and on $X//\T$.
See \S 5.1 of  \cite{BCK2} for more details on $X//\G$, $X//\T$,
and the $\bS$-actions on them. As before, we let
$\CC[\lambda]=\CC[\lambda_1,\dots,\lambda_n]=H^*(B\bS,\CC)$, with quotient field
$\CC(\lambda)$.
Our goal in this section is to prove

\begin{Thm} \label{mainthm} Conjecture \ref{conj} holds for

$(a)$ the usual Gromov-Witten theory of $Fl(k_1,...,k_r,n)$.

$(b)$ the $\bS$-equivariant Gromov-Witten theory of $Fl(k_1,...,k_r,n)$.
\end{Thm}

\begin{Rmk} Note that part $(a)$ follows from $(b)$ by taking
the non-equivariant limit $\lambda_1=\dots =\lambda_n=0$ of the potential
functions.
\end{Rmk}

Our strategy for proving Theorem \ref{mainthm} is to use
reconstruction theorems to reduce the statement to a comparison for $1$-point
invariants which was established in \cite{BCK1}, \cite{BCK2}. Typically, reconstruction
results for Gromov-Witten invariants work under the assumption that the cohomology ring is generated
by divisors. Our observation here is that in the
torus-equivariant setting, this assumption needs only to hold after localization.
This enlarges the class of varieties for which reconstruction is applicable.
We begin with a simple lemma.

\begin{Lemma}\label{loc} Let $\PP ^N$ be acted by a torus $\bS$
and let $Y$ be an $\bS$-invariant smooth subvariety. Suppose that the
natural map $H^*((\PP ^N)^\bS)\ra H^*(Y^\bS)$ is surjective (for
example, this is true when the $\bS$-fixed locus $(\PP^N)^\bS$ is isolated).
Then the localized equivariant cohomology ring $H^*_\bS(Y,\CC)\otimes_ {\CC[\lambda ]}
\CC(\lambda)$ is generated (as a $\CC(\lambda)$-algebra) by divisors, i.e., by
$\{c_1(L)\otimes 1\;|\; L\in\Pic^\bS(Y)\}$.
\end{Lemma}

\begin{proof}  There is a commutative diagram
\[ \begin{CD} H^*_\bS(\PP ^N ) @<<< H^*((\PP ^N)^\bS)\otimes\CC[\lambda] \\
@VVV @VVV \\
H^*_\bS(Y) @<<< H^*(Y^\bS)\otimes\CC[\lambda]
\end{CD}\]
After tensor product with $\CC (\lambda)$, the horizontal maps are
isomorphism by the localization theorem. \end{proof}

\medskip

It is well-known that $X//\G=Fl(k_1,\dots ,k_r,n)$ admits
an $\bS$-equivariant embedding into a product of projective spaces on which $\bS$ acts
with isolated fixed points. By Lemma \ref{loc}, the localized equivariant cohomology
$$H^*_{\bS}( X//\G,\CC)\otimes_{\CC[\lambda]}\CC(\lambda)$$
is generated by divisor classes. Note that this is false in general without localization.
For example the equivariant cohomology of Grassmannians is {\em not} generated by divisors.
(On the other hand, since $X//\T$ is a toric variety, both the usual and $\bS$-equivariant
cohomology rings are already generated by divisors.)

Before going into the details of the proof, it is useful to discuss the base-change of
Novikov rings $(3.2.1)$ in the particular case of flag manifolds. By choosing the
usual Schubert basis in $H_2(Fl(k_1,\dots,k_r,n),\ZZ)$, the semigroup of effective curve classes
is identified with $(\ZZ_+)^r$. We write $d=(d_1,\dots,d_r)$ for the typical element in
this semigroup. The Novikov ring is simply the power series ring $\CC[[Q_1,\dots,Q_r]]$.
Similarly, effective curve classes on the toric variety $X//\T$ are described by tuples
of non-negative integers
$$\tilde{d}=(d_{11},\dots, d_{1k_1},\dots,d_{r1},\dots,d_{rk_r})$$
and the Novikov ring is identified with $\CC[[Q_{ij}|1\leq i\leq r,\; 1\leq j\leq k_i]]$.
A class $\tilde{d}$ is a lift of $d$ if and only if
$$d_i=\sum_{j=1}^{k_i}d_{ij},\;\;\; i=1,\dots, r.$$
Finally, $\epsilon(d)=\sum_{i=1}^r(k_i-1)d_i\; (\mathrm{mod}2)$.
Hence the projection $(3.2.1)$ of Novikov rings is
\begin{equation}\label{flag specialization}
p:\CC[[Q_{ij}|1\leq i\leq r,\; 1\leq j\leq k_i]]\longrightarrow
\CC[[Q_1,\dots,Q_r]],\;\;\;\; Q_{ij}\mapsto (-1)^{(k_i-1)}Q_i.
\end{equation}

\subsection{Kontsevich-Manin reconstruction and
reduction to 2-point invariants} This step involves an equivariant version
of the Kontsevich-Manin Reconstruction Theorem.
In its original formulation
\cite{KM}, the Reconstruction Theorem states that if the
cohomology ring $H^*(Y,\CC)$ is generated by divisors, then all Gromov-Witten invariants
of $Y$ can be reconstructed from $3$-point invariants for which at least one
insertion is a divisor class. These in turn are expressed in
terms of $2$-point invariants by using the divisor equation in Gromov-Witten
theory. We give here an extension of reconstruction to the $\bS$-equivariant setting.

\begin{Lemma} \label{eqKM} Let $Y$ be a smooth complex projective variety with $\bS$-action.
Let $$P:=Spf((N(Y)\otimes\CC[\lambda])[[H^*_\bS(Y,\CC)^\vee]]),$$
endowed with the formal $\bS$-equivariant Frobenius structure
$(P,\star,g,1,\fE_{\bS})$ defined by the equivariant Gromov-Witten
potential $F^Y$. Let $t=(t_0,t_1,\dots,t_r,\dots,t_{m-1})$ be the
flat coordinates defined by a basis of $H^*_\bS(Y,\CC)$, such that
$t_{\mathrm{small}}=(t_1,\dots,t_r)$ are the coordinates on the
small parameter space $H^2_\bS(Y,\CC)$. Let $G\in\cO_P$ be another
formal function satisfying the WDVV equations. Assume that:

\noindent $(i)$ In flat coordinates $G$ has the form $(\ref{qc-type})$
$$G=G_{cl}+\sum_{\beta\in E,\beta\neq 0}
Q^{\beta}e^{\beta\cdot t_{\mathrm{small}}}G_{\beta},
$$
with $G_{\beta}\in\CC[\lambda][[t_{r+1},\dots,t_{m-1}]]$ and $G_{cl}$ a
cubic polynomial in the $t_i$'s (with coefficients in
$\CC[\lambda]$).

\noindent $(ii)$ $\cL_{\fE_\bS}(G)=(3-\dim(Y))G$.

\noindent $(iii)$ $G_{cl}=F^Y_{cl}$.

\noindent $(iv)$ $\left.\partial_{t_i}\partial_{t_j}G\right|_{t_{\mathrm{small}}}=
\left.\partial_{t_i}\partial_{t_j}F^Y\right|_{t_{\mathrm{small}}}$, for all $i,j$.

\noindent $(v)$ The localized equivariant
cohomology ring $H^*_\bS(Y,\CC)\otimes_{\CC[\lambda]}\CC(\lambda)$ is generated
by $H^2_\bS(Y,\CC)$ as a $\CC(\lambda)$-algebra.

Then $G=F^Y$.
\end{Lemma}

\begin{proof} Let $P_{(\lambda)}$ be the
$\bS$-equivariant Frobenius manifold defined by the localized
Gromov-Witten theory of $Y$ (see \S 2.1). The function $G$ defines
a formal $\bS$-equivariant Frobenius structure $(P,\circ,g,1,
\fE_{\bS})$ over $N(Y)\otimes\CC[\lambda]$, and a localized
Frobenius structure over $N(Y)\otimes \CC(\lambda)$ as well, by
viewing it as a formal function on $P_{(\lambda)}$ via the natural
(injective!) localization map $\iota:\cO_P\ra\cO_{P_{(\lambda)}}$.
It suffices to check that the localized potentials
$F^Y_{(\lambda)}=\iota(F^Y)$ and $G_{(\lambda)}=\iota(G)$ are
equal. The assumptions $(i)$-$(iii)$ hold for the localized
potentials as well (where in $(i)$ we replace $\CC[\lambda]$ by
$\CC(\lambda)$).

In the conformal case, a formal Frobenius structure satisfying
$(i)$ and $(ii)$ is said to be of {\em qc-type} in \cite{Manin}.
Such structure has \lq\lq cup product", defined by
$$G_{cl}=\frac{1}{6}g((\sum t_i\partial_{t_i})\cup (\sum
t_i\partial_{t_i}),\sum t_i\partial_{t_i} )$$ and \lq\lq
correlators"
$$I_{0,n,\beta}(\partial_{t_{i_1}},\dots , \partial_{t_{i_n}})=
\partial_{t_{i_1}}\dots  \partial_{t_{i_n}}G_\beta|_{t=0}$$
which satisfy the analogue of the divisor axiom in Gromov-Witten theory. See \cite{Manin}, \S 5.4.
The same will hold for the Frobenius structure defined by our potential $G$,
or for its localized version.
We may call them Frobenius structures of equivariant qc-type.

Assumption $(v)$ and the usual Kontsevich-Manin reconstruction imply that the localized GW-potential
$F^Y_{(\lambda)}$ is determined recursively by
$\partial_{t_i}\partial_{t_j}F^Y_{(\lambda)}|_{t_{\mathrm{small}}}$.
The proof only uses properties of Gromov-Witten invariants which are shared by
the correlators of
any Frobenius structure of equivariant qc-type, hence it will work in the abstract case as well.

Assumption $(iii)$ says that the abstract cup product coincides
with the usual one on cohomology. By the above discussion
reconstruction applies and we find that $G_{(\lambda)}$ is
determined recursively by
$\partial_{t_i}\partial_{t_j}G_{(\lambda)}|_{t_{\mathrm{small}}}$,
{\em with the same recursion coefficients} in $\CC(\lambda)$ as
those for $F^Y_{(\lambda)}$. By assumption $(iv)$, we are done.
\end{proof}

\bigskip

We go back now to the proof of Theorem 4.1.
We intend to apply
Lemma \ref{eqKM} to $Y=Fl$, $G=\varphi^*F'$.
Note that assumption $(v)$ holds by Lemma \ref{loc},
assumption $(i)$ holds by (\ref{qc-type'}), while
$(ii)$ and $(iii)$ are immediate from the construction of $F'$ in $\S 3.2 - \S 3.5.$
Hence the Theorem will be proved if we can show that $(iv)$ holds as well, i.e.,
\begin{equation}\label{restricted partials equal}
\left.\partial_{t_i}\partial_{t_j}F^{Fl,\bS}\right|_{t_0=t_{r+1}=\dots=t_{m-1}=0}
=\varphi ^*\left(
\left.\partial_{s_i}\partial_{s_j}F'\right|_{s_0=s_{r+1}=\dots=s_{m-1}=0}\right).
\end{equation}
Recall that (in the notation of \S 2.2)
$$\partial_{t_i}\partial_{t_j}F^{Fl,\bS}(t)=\langle\langle \sigma_i,\sigma_j\rangle\rangle .$$
Setting $t_0=t_{r+1}=\dots=t_{m-1}=0$ and using the divisor axiom we get
\begin{equation}\label{conversion1}
\left.\partial_{t_i}\partial_{t_j}F^{Fl,\bS}\right|_{t_0=t_{r+1}=\dots=t_{m-1}=0} =
\sum_{d=(d_1,\dots,d_r)}\prod_{l=1}^r (Q_le^{t_l})^{d_l}\langle \sigma_i,\sigma_j\rangle^{Fl,\bS}_{0,2,d}.
\end{equation}
On the other hand, since
$$\partial_{s_k}\partial_{s_i}\partial_{s_j}F'=\partial_{s_k}\left(\left.\left(\partial_{\tilde{t}_i\cup\omega}
\partial_{\tilde{t}_j\cup\omega}F\right)\right|_N\right)$$
it follows that
\begin{equation}\label{partials equal}
\partial_{s_i}\partial_{s_j}F'=\left.\left(\partial_{\tilde{t}_i\cup\omega}
\partial_{\tilde{t}_j\cup\omega}F\right)\right|_N
\end{equation}
up to a constant (in the base ring). By adding appropriate quadratic terms to $F'$,
we may assume that (\ref{partials equal})
holds exactly. (Recall that $F$ is the $\bS$-equivariant Gromov-Witten potential
of $X//\T$ with the Novikov variables specialized as in (\ref{flag specialization}).)
Moreover, the first Chern class of the flag manifold satisfies
\begin{equation*}
\int_{d}c_1(T(Fl))=\sum_{l=1}^rd_l(k_{l+1}-k_{l-1})\geq \min_l\{k_{l+1}-k_{l-1}\}
\geq 2.
\end{equation*}
Therefore the specialization of the left-hand side of
(\ref{partials equal}) to
$s_0=s_{r+1}=\dots=s_{m-1}=0$ is equal to its specialization at
$\tilde{t}_{0}=\tilde{t}_{r+1}=
\dots=\tilde{t}_{m-1}=0$ by Lemma \ref{small space} $(i)$.
Using the divisor axiom as above in the right hand side of (\ref{partials equal}),
we conclude that
%\begin{equation}
\begin{multline}\label{conversion2}
\left.\partial_{s_i}\partial_{s_j}F'\right|_{s_{r+1}=\dots=s_{m-1}=0}=\\
\sum_{d=(d_1,\dots,d_r)}\prod_{l=1}^r (Q_l e^{\tilde{t}_l})^{d_l}
\left(\sum_{\tilde{d}\mapsto d}(-1)^{\sum(k_l-1)d_l}
\langle \gamma_i\cup\omega,\gamma_j\cup\omega\rangle^{X//\T,\bS}_{0,2,\tilde{d}}\right).
\end{multline}
%\end{equation}
Using (\ref{conversion1}) and (\ref{conversion2}), the proof of
(\ref{restricted partials equal}), and therefore of
Theorem \ref{mainthm},
is reduced to checking the following identity among $2$-point invariants:

\begin{equation}\label{2 pt primary}
\langle \sigma_i,\sigma_j\rangle^{Fl,\bS}_{0,2,d}=\sum_{\tilde{d}\mapsto d}(-1)^{\sum(k_l-1)d_l}
\langle \gamma_i\cup\omega,\gamma_j\cup\omega\rangle^{X//\T,\bS}_{0,2,\tilde{d}}.
\end{equation}

\subsection{Lee-Pandharipande reconstruction and reduction to $1$-point invariants}\label{LP reconstruction}
There is another reconstruction theorem, due to Lee and
Pandharipande \cite{LP1}, and independently to Bertram and Kley,
which reduces in certain cases computations of (descendant)
$GW$-invariants with any number of insertions to 1-point
descendants. In fact, Lee and Pandharipande deduce the
reconstruction result from universal relations they found among
divisors in the Picard group of the moduli space
$\overline{M}_{0,2}(\PP^N,d[line])$ of $2$-pointed stable maps to
$\PP^N$. We establish first a straightforward equivariant
extension of their divisor relation.

Let $Y$ be a projective variety with $\bS$-action. Let $\overline{M}_{0,2}(Y,\beta)$
be the moduli space of 2-pointed genus zero stable maps with evaluation maps
$$ev_1,ev_2:\overline{M}_{0,2}(Y,\beta)\lra Y.$$
The moduli space inherits an $\bS$-action and the evaluation maps are equivariant.
Let $\psi=\psi_1$ be the $\bS$-equivariant
first Chern class of the line bundle on $\overline{M}_{0,2}(Y,\beta)$
with fiber $T^*_{x_1}C$ over the point $[f:(C,x_1,x_2)\ra Y]$.

There is a \lq\lq boundary divisor" $D_{1,\beta_1|2,\beta_2}$ in
$\overline{M}_{0,2}(Y,\beta)$ corresponding to maps with reducible
domains and splitting type $\beta_1+\beta_2=\beta$. It is obtained
as the image of the ($\bS$-equivariant) gluing morphism
$$j_{\beta_1,\beta_2}:\overline{M}_{0,\{x_1,\bullet\}}(Y,\beta_1)\times_Y\overline{M}_{0,\{x_2,\bullet\}}
(Y,\beta_2)\lra \overline{M}_{0,2}(Y,\beta)$$ and its virtual
fundamental class in the equivariant Chow group
$A_*^\bS(\overline{M}_{0,2}(Y,\beta),\QQ)$ is defined as the
push-forward of
$$[\overline{M}_{0,\{x_1,\bullet\}}(Y,\beta_1)]^{vir}
\boxtimes [\overline{M}_{0,\{x_2,\bullet\}}(Y,\beta_2)]^{vir}.
$$

\begin{Lemma} \label{eqLP} For all $L\in\Pic^{\bS}(Y)$, the relation
\begin{multline}\notag
ev_2^*(L)\cap[\overline{M}_{0,2}(Y,\beta)]^{vir}-
\left(ev_1^*(L)+(\int_\beta L)\psi\right)\cap[\overline{M}_{0,2}(Y,\beta)]^{vir}
+\\
\sum_{\beta_1+\beta_2=\beta}(\int_{\beta_2} L)\cap[D_{1,\beta_1|2,\beta_2} ]^{vir}=0
\end{multline}
holds in $A_*^\bS(\overline{M}_{0,2}(Y,\beta),\QQ )$.
\end{Lemma}

\begin{proof} As in \cite{LP1}, since the relation is linear in $L$ and
the equivariant Picard group is spanned over $\QQ$ by $\bS$-equivariant very ample line
bundles, the Lemma will follow from the case $Y=\PP^N$, $\beta=d[line]$
and the stronger statement
\begin{equation}\label{LPrel}
ev_2^*(L)-
ev_1^*(L)-(\int_\beta L)\psi
+\sum_{\beta_1+\beta_2=\beta}(\int_{\beta_2} L)D_{1,\beta_1|2,\beta_2} =0
\end{equation}
in $\Pic^\bS(\PP^N)$.

The relation (\ref{LPrel}) holds after passing to the
non-equivariant limit $\lambda_i=0$ by Theorem 1 in \cite{LP1}.
Therefore the left-hand side is a linear polynomial in the
$\lambda_i$'s and the corresponding equivariant line bundle is
just a trivial bundle twisted by a character of $\bS$. To check
that this character is trivial, it suffices to restrict to {\em
any} $\bS$-fixed point of $\overline{M}_{0,2}(\PP^N,d[line])$.
There are many possible choices of fixed points that will work.
One particular such for which the computation is very easy is the
point corresponding to a stable map with domain $C\cup D$ (the
union of two irreducible components) such that $x_1,x_2\in C$ and
$f:C\cup D\ra\PP^N$ collapses $C$ to a fixed point $p\in\PP^N$ and
maps $D$ with degree $d$ onto an $\bS$-invariant line in $\PP^N$
joining $p$ to another fixed point $q$, such that the map is
totally ramified at $q$. The classes $\psi$ and
$D_{1,\beta_1|2,\beta_2}$ vanish when restricted to this point,
while $ev_1^*L$ and $ev_2^*L$ have the same restriction. Relation
(\ref{LPrel}) and hence the Lemma are proved. \end{proof}

\medskip

We will use Lemma \ref{eqLP} to obtain a reconstruction result in
the context of the abelian/nonabelian correspondence. Recall that
descendant (genus $0$) Gromov-Witten invariants of a smooth
projective $Y$ are defined by
$$\langle \tau_{a_1}(\gamma_1),\dots,\tau_{a_n}(\gamma_n)\rangle^Y_{0,n,\beta}:=
\int_{[M_{0,n}(Y,{\beta})]^{vir}}\prod_i\psi_i^{a_i}ev_i^*(\gamma_i),
$$
where $\gamma_i\in H^*(Y)$ and $\psi_i$ are the first Chern
classes of the cotangent line bundles at the marked points. The
definition extends to torus-equivariant descendants (which will be
$\CC[\lambda]$-valued). We establish first an auxiliary vanishing
result for certain descendant invariants of $X//\T$.

Let $X//\G$, $X//\T$, $\bS$ etc. be as in the setting \S 3.1. Let $\beta\in H_2(X//\G,\ZZ)$
be fixed. Consider the moduli space
$$\mathcal{M}_{\beta}:=\coprod_{\tilde{\beta}\mapsto\beta}\overline{M}_{0,n}(X//\T,\tilde{\beta})
$$
with the obvious evaluation maps $ev_i:\cM_\beta\lra X//\T$,
$i=1,\dots, n$ and virtual class $[\cM_{\beta}]^{vir}$. Note that
$$H^*_{\bS}(\cM_{\beta},\CC)\cong\bigoplus_{\tilde{\beta}\mapsto\beta}
H^*_\bS(\overline{M}_{0,n}(X//\T,\tilde{\beta}),\CC).$$ Introduce
\lq\lq psi-classes" on $\cM_{\beta}$ by
$$\psi_i:=\sum_{\tilde{\beta}\mapsto\beta}\psi_{i,\tilde{\beta}}$$
and define for cohomology classes $\gamma_1,\dots,\gamma_n\in
H^*_\bS(X//\T)$.
\begin{equation*}\begin{split}
I_{n,\beta}(\tau_{a_1}(\gamma_1),\dots,\tau_{a_n}(\gamma_n))&:=
(-1)^{\epsilon(\beta)}
\int_{[\cM_{\beta}]^{vir}}\prod_i\psi_i^{a_i}ev_i^*(\gamma_i)\\
&=(-1)^{\epsilon(\beta)}
\sum_{\tilde{\beta}\mapsto\beta}
\langle \tau_{a_1}(\gamma_1),\dots,\tau_{a_n}(\gamma_n)
\rangle_{0,n,\tilde{\beta}}^{X//\T,\bS}.
\end{split}\end{equation*}
Recall that the intersection form is non-degenerate on
the$\bW$-anti-invariant subspace $H^*_\bS(X//\T)^a$. We denote the orthogonal
complement by $\left(H^*_\bS(X//\T)^a\right)^{\perp}$.

\begin {Lemma} \label{vanishing}
If $\tilde{\sigma}_1,\dots,\tilde{\sigma}_{n-1}$ are $\bW$-invariant lifts
of classes $\sigma_i$ in $H^*_\bS(X//\G)$ and $\gamma\in\left(H^*_\bS(X//\T)^a\right)^{\perp}$,
then
$$I_{n,\beta}(\tau_{a_1}(\tilde{\sigma}_1\cup\omega),\tau_{a_2}(\tilde{\sigma}_2),\dots,
\tau_{a_{n-1}}(\tilde{\sigma}_{n-1}),\tau_{a_n}(\gamma))=0.$$
\end{Lemma}

\begin{proof} The $\bW$-action on $X//\T$ induces a $\bW$-action on $\cM_{\beta}$, by
composing stable maps with the automorphisms in $\bW$. The evaluation maps are easily seen to be
$\bW$-equivariant. Note also that the psi-classes are $\bW$-invariant. Hence the class
$$(ev_n)_*
\left(ev_1^*(\omega)\prod_{i=1}^{n}\psi_i^{a_i}\prod_{i=1}^{n-1}ev_i^*(\tilde{\sigma}_i)
\cap[\cM_\beta]^{vir} \right)$$ is $\bW$-anti-invariant. The Lemma
now follows from the projection formula. \end{proof}

\medskip

\begin{Prop} Let $X//\G$, $X//\T$, $\bS$ be as in the setting \S 3.1.
Assume that the localized equivariant cohomology $H^*_\bS(X//\G,\CC)\otimes_{\CC[\lambda]}\CC(\lambda)$
is generated as a $\CC(\lambda)$-algebra by
divisors (that is, by $c_1(L)\otimes 1$ for $L\in\Pic^\bS(X//\G)$). Let $\sigma_i,\sigma_j$ be any equivariant
cohomology classes on $X//\G$, with $\bW$-invariant lifts $\tilde{\sigma}_i,\tilde{\sigma_j}$ to $X//\T$. If
the identity
\begin{equation}\label{2 pt descendants}
\langle \tau_a(\sigma_i),\sigma_j\rangle^{X//\G,\bS}_{0,2,\beta}=
\sum_{\tilde{\beta}\mapsto \beta}(-1)^{\epsilon(\beta)}
\langle \tau_a(\tilde{\sigma}_i\cup\omega),\tilde{\sigma}_j\cup\omega\rangle^{X//\T,\bS}_{0,2,\tilde{\beta}}
\end{equation}
holds for $\sigma_j=1$, then it holds in general.
\end{Prop}

\begin{proof} It is enough to prove the Proposition for a {\em fixed} choice
of lifts of cohomology classes on $X//\G$ to $X//\T$.

It follows immediately from Martin's integration formula (3.1.3) that
\begin{equation}\label{lifting}
\widetilde{\sigma'\cup\sigma''}\cup\omega=\widetilde{\sigma'}\cup\widetilde{\sigma''}\cup\omega
\end{equation}
for any $\sigma',\sigma''\in H^*_\bS(X//\G,\CC)$ (see e.g. Cor. 2.3
in \cite{BCK1} for an argument).

Assume first that the equivariant cohomology ring of $X//\G$ is
generated by divisors without localization (this happens for
example when $X//\G$ is the complete flag manifold
$Fl(1,2,\dots,n-1,n)$).  Using Lemma \ref{eqLP} and the splitting
axiom for GW-invariants we find that $\langle
\tau_a(\sigma_i),\sigma_j\rangle^{X//\G,\bS}_{0,2,\beta}$ is
expressed recursively (with $\CC[\lambda]$-coefficients) in terms
of invariants $\langle
\tau_{a'}(\sigma'),1\rangle^{X//\G,\bS}_{0,2,\beta'}$ (these can
be further reduced to 1-point descendants by the fundamental class
axiom for GW-invariants). This is just the reconstruction of
Lee-Pandharipande.

Recall the notation
$I_{2,\beta}(\tau_a(\tilde{\sigma}_i\cup\omega),\tilde{\sigma}_j\cup\omega)$
introduced above for the right-hand side of the identity (\ref{2
pt descendants}). The divisor relation in Lemma \ref{eqLP} can be
extended in an obvious manner to the moduli space $\cM_{\beta}$
for $\bW$-invariant lifts $\tilde{L}$ of line bundles
$L\in\Pic^\bS(X//\G)$. By Lemma \ref{vanishing}, the
reconstruction procedure applies to the invariants
$I_{2,\beta}(\tau_a(\tilde{\sigma}_i\cup\omega),\tilde{\sigma}_j\cup\omega)$
and (by (\ref{lifting}) and the equality
$\epsilon(\beta_1+\beta_2)=\epsilon(\beta_1)+\epsilon(\beta_2)$)
it expresses them in terms of
$I_{2,\beta'}(\tau_{a'}(\tilde{\sigma'}\cup\omega),\omega)$ with
the {\underline {same}} recursion coefficients. The Proposition is
proved in this case.

In the general case the same argument will work word for word,
except that the recursion coefficients will now be rational
functions rather than polynomials in the $\lambda_i$'s.
\end{proof}

\begin{Rmk} In view of Lemma \ref{vanishing}, one might be tempted to try to
extend the  version of Lee-Pandharipande reconstruction above to
descendants with any number of insertions and $\psi$-classes at
all points. However, this is not possible, because an analogue of
the fundamental class axiom does not hold for the invariants
$I_{n,\beta}$ (indeed, in general
$I_{n,\beta}(\tilde{\sigma}_1\cup\omega,\tilde{\sigma}_2,\dots,\tilde{\sigma}_{n-1},\omega)\neq
0$). This is the reason for which the \lq\lq twisting" by
$\star\omega$ is necessary.

\end{Rmk}

\begin{Cor} The following identity holds between Gromov-Witten invariants of $X//\G=Fl(k_1,\dots,k_r,n)$
and those of the abelian quotient $X//\T$: for any $d=(d_1,\dots,d_r)\in H_2(Fl,\ZZ)$, any $a\geq 0$
and any equivariant cohomology classes $\sigma_i,\sigma_j$ on $Fl$, with lifts $\gamma_i,\gamma_j$ respectively,

$$\langle \tau_a(\sigma_i),\sigma_j\rangle^{Fl,\bS}_{0,2,d}=\sum_{\tilde{d}\mapsto d}(-1)^{\sum(k_a-1)d_a}
\langle \tau_a(\gamma_i\cup\omega),\gamma_j\cup\omega\rangle^{X//\T,\bS}_{0,2,\tilde{d}}$$
\end{Cor}

\begin{proof} By Lemmas \ref{loc} and \ref{eqLP}, it suffices to check that
\begin{equation}\label{1 pt}
\langle \tau_a(\sigma_i),1\rangle^{Fl,\bS}_{0,2,d}=\sum_{\tilde{d}\mapsto d}(-1)^{\sum(k_a-1)d_a}
\langle \tau_a(\gamma_i\cup\omega),\omega\rangle^{X//\T,\bS}_{0,2,\tilde{d}}.
\end{equation}

This is (essentially) proved in \cite{BCK1}, \cite{BCK2}. However, since the actual statement
is explicitly written  (see formula (5) on p. 124 and Remark on p. 125 in \cite{BCK1})
only for Grassmannians and non-equivariant
invariants, we should say a few words here.

For the general flag manifold, a correspondence between the
equivariant \lq\lq small" $J$-functions of $Fl$ and $X//\T$ is
given by Theorem 1 in \cite{BCK2} (see the next section below for
more about $J$-functions). Reading the argument on p. 124-125 in
\cite{BCK1} backwards\footnote{The specialization of the
$t_i$-variables there corresponds exactly to our specialization
(\ref{flag specialization}) of the Novikov variables $Q_i$ here.},
the equality (\ref{1 pt}) follows from the $J$-functions
correspondence, provided that for any factorization
$$\omega=\left(\sqrt{\frac{(-1)^{|\Phi_+|}}{|\bW|}}\right)(\prod_{\alpha\in A}c_1^\bS(L_{\alpha}))\cup
(\prod_{\alpha\in \Phi_+\setminus A}c_1^\bS(L_\alpha)) $$
we have
\begin{equation}\label{quantum omega = omega}
\left(\sqrt{\frac{(-1)^{|\Phi_+|}}{|\bW|}}\right)
(\prod_{\alpha\in A}c_1^\bS(L_{\alpha}))\star_{\mathrm{small}}
(\prod_{\alpha\in \Phi_+\setminus A}c_1^\bS(L_\alpha))=\omega,
\end{equation}
where $\star_{\mathrm{small}} $ is the small equivariant quantum
product on $X//\T$, restricted to $H^2_\bS(X//\T,\CC)^\bW$ and
with the Novikov variables specialized as in (\ref{flag
specialization}). By a simple degree counting, this last equality
is always true when $X//\G$ (and hence $X//\T$, cf. \S 3.4) is a
Fano variety. Indeed, the left-hand side of (\ref{quantum
omega = omega}) is $\bW$-anti-invariant, homogeneous, and
of the form
$$\omega +\;\mathrm{quantum\; corrections}.$$ However, $\omega$ is the
unique class of lowest degree in $H^*_\bS(X//\T,\CC)^a$, and in
the Fano case the quantum parameters have positive degree. Hence
the quantum corrections must vanish. \end{proof}

\medskip

It remains to observe that relation (\ref{2 pt primary}) is a
special case of the Corollary to conclude the proof of Theorem
\ref{mainthm}. \qed

\bigskip

Note that the only instance in this section where we have used that $X//\G$ is
a flag manifold was in quoting the small $J$-function correspondence from
\cite{BCK2}. In other words, we have proved

\begin{Thm} Let $X,\G,\T,\bS$ etc. be as in the setting \S 3.1. Assume that
$X//\G$ is Fano of index $\geq 2$ and that its equivariant cohomology is
generated by divisors after localization. Then Conjecture \ref{conj} holds
if and only if equation (\ref{1 pt}) holds, if and only if the abelian/nonabelian
correspondence for small $J$-functions holds.
\end{Thm}

\noindent A similar statement holds if we only assume that $c_1(T(X//\G))$ is nef,
by using Lemma \ref{small space} $(ii)$ in the argument just above
equation (\ref{conversion2}).
However, the change of coordinates $s(\tilde{t})$ will be nontrivial even
for the restriction to subspace $\{s_{r+1}=\dots=s_{m-1}=0\}$,
and coincides with the change of coordinates in the abelian/nonabelian
correspondence for small $J$-functions (see Conjecture 4.3 in \cite{BCK2}).
This is precisely analogous to the mirror theorem \cite{Giv1} for hypersurfaces
in projective space.
We leave the precise formulation for the interested reader.

\section{The abelian/nonabelian correspondence for
$J$-functions}\label{J-function}

Our goal in this section is
to explain why Conjecture \ref{conj} is equivalent to (an extension to the big
parameter space of) the correspondence between
the $J$-functions of $X//\G$ and $X//\T$ proposed in \cite{BCK2}, Conjecture 4.3.
In particular,
by Theorem \ref{mainthm} and Corollary \ref{J-I correspondence} below,
the correspondence holds for the flag manifolds
$Fl(k_1,\dots,k_r,n)$.

\subsection{Deformed flat coordinates}
First we recall the definition of deformed flat coordinates
following Dubrovin \cite{Du1, Du2, Du3}. Let $M$ be a Frobenius
manifold (say, analytic, for simplicity), with Euler vector field.
There is a deformed flat connection $\nabla ^z$ on $\Theta_M$
given by
$${\nabla}^z
_\xi \eta :=\nabla _\xi \eta-z^{-1}\xi \star\eta$$ (see p. 189 and
p. 323 of \cite{Du1} and also \cite{Du3}; however, we follow
Givental for the convention on $z$). By identifying the cotangent
sheaf $\Omega^1_M$ and the tangent sheaf $\Theta_M$ via the flat
metric, a deformed flat connection is induced on $\Omega^1_M$. A
coordinate system $J_i$ of $M$ is called a {\em deformed flat
coordinate system} if $dJ_i$ are horizontal sections. In other
words, $J_i$ form a complete solution space to the second order
linear PDE system
\begin{eqnarray}\label{J-eq} z\partial _{t_i}
\partial _{t_j} J = \sum _\gamma c^k _{ij} \partial _{t_k} J
\end{eqnarray}
where $t_i$ are flat coordinates and $c_{ij}^k$ are structure
constants of multiplications, i.e., $\partial _{t_i}\star\partial
_{t_j} =\sum _k c^k_{ij} \partial _{t_k}$.

Suppose that the potential
function $F$ (defined up to quadratic terms) for the Frobenius
structure is of the form $F=F_c+F_q$, with $F_c$ a cubic form of
the flat coordinates $t_i$ and $F_q \in
\CC[[q_1,...,q_r,t_{r+1},...,t_R]]$ such that $q_i=e^{t_i}$ and
$F_q\equiv 0$ modulo the ideal $(q_1,...,q_r)$. (cf. \ref{qc-type})

Consider the normalization
condition
$$\sum J_{i}\partial _{t_i}
\equiv ze^{\bt/z} = z\partial _{t_0}+ \bt+
O(z^{-1})\;(\mathrm{mod}(q_1,\dots,q_r)),$$ where $\bt=\sum t_i
\partial _{t_i}$, the products of vector fields in the exponential
are the \lq\lq cup" products (determined by $F_{cl}=(1/6)g(\bt\cup
\bt,\bt)$) and $1=\partial _{t_0}$. The normalization uniquely
determines deformed flat coordinates once the flat coordinates are
chosen (see Lemma 2.2 of \cite{Du2}). We will call $\sum J_i
\partial _i $ the {\em $J$-function} if it is normalized as above.

\subsection{$J$-functions in quantum cohomology}
The $J$-function plays an important role in Gromov-Witten Theory.
Let $Y$ be a projective algebraic manifold. Then the $J$-function
$J_Y$ for the (formal) Frobenius structure defined by the quantum
cohomology of $Y$ can be constructed using descendant
Gromov-Witten invariants. Let $\{\phi_i\}$ be a homogeneous basis
of $H^*(Y)$, with Poincar\'e dual basis $\{\phi^i\}$. Let
${\bt}:=\sum_it_i\phi_i$. $J_Y$ coincides with the assignment

\begin{equation}\label{Jfcn}
H^*(Y)\ni {\bt} \mapsto z+\bt+ \sum _{i} \phi ^{i}
\langle\langle\frac{\phi _i}{z-\psi}\rangle\rangle \ \ \in z +
\bt+ \mathcal{H}_-
\end{equation}
cf. \cite{CG, SGF}, where $\mathcal{H}_-=\frac{1}{z}H^*(Y)\otimes
_{\CC} N[Y] [[\frac{1}{z}]]$. Precisely speaking, $J_Y$ is an
element in $R[[K^\vee ]]((1/z))$.

Here we use the double-bracket notation introduced in \S 2.2, so
that
\begin{equation*}
\langle\langle\frac{\phi _i}{z-\psi}\rangle\rangle= \sum_{\beta\in
NE_1}Q^\beta\sum_{n\geq
0}\frac{1}{n!}\int_{[\overline{M}_{0,n+1}(Y,\beta)]^{vir}}
\frac{ev_1^*(\phi_i)}{z-\psi} ev_2^*(\bt)\dots ev_{n+1}^*(\bt)
%=\sum_{\beta}Q^{\beta}J_{i,\beta}
\end{equation*}
where $\psi=\psi_1$ and $1/(z-\psi)$ is formally expanded as a
geometric series.

The normalization condition $$J_Y(\bt ,z)\equiv ze^{\bt/z}
$$
modulo quantum corrections follows from
the well-known result
$$\int
_{\overline{M}_{0,n}}\psi _1^{l_1}...\psi _n^{l_n} =
(n-3)!/l_1!...l_n!\;\;\;\;\; {\mathrm{if}} \;\;\sum l_i = n-3.$$
(Note that in the paper \cite{BCK2} $J_Y(t,z)/z$ is used for
$J$-function, i.e., a different normalization.)

\subsection{The abelian/nonabelian correspondence}
Let $X,\G,\T$ be as in the setting \S 3.1. (For simplicity, we do not
consider the equivariant theory here; the interested reader can readily make the
necessary modifications to cover this case as well.) We have the
$\bW$-induced Frobenius structure over the Novikov ring $N(X//\G)$
constructed in \S3.2 - \S 3.6. We will keep the notations, and make liberal
use of all its properties proved there. Moreover,
{\em from now on, we assume
that Conjecture \ref{conj} holds for $X//\G$ and $X//\T$.}

\bigskip

If
$J_{X//\G}=\sum_{i=0}^{m-1}J_{i,X//\G}(t_0,\dots,t_{m-1},z)\sigma_i$
is the $J$ function of $X//\G$, as given by (\ref{Jfcn}), put
$$\tilde{J}_{X//\G}(t,z):= \sum_{i=0}^{m-1}J_{i,X//\G}(t_0,\dots,t_{m-1},z)\gamma_i.$$
(Recall that $\gamma_i$'s are chosen $\bW$-invariant lifts of the $\sigma_i$'s.)

\begin{Lemma}\label{d omega}
$\tilde{J}_{X//G} (t,z)\cup \omega =\left. (z\partial _\omega
J_{X//T})\right|_{Q^{\tilde{\beta}
}=(-1)^{\epsilon(\beta)}Q^\beta, N} (\varphi (t),z).$
\end{Lemma}

\begin{proof} Both sides satisfy the normalization
condition $J\equiv ze^{t/z}\cup \omega$ modulo quantum
corrections. Therefore it suffices to check that $\{\partial
_\omega J_i \}_i$ forms a deformed flat coordinate system for $(N,
\circ, ^\omega g, e, \mathfrak{E})$ if $J_\delta$ is a deformed
flat coordinate for $(M,\star,g,e,\mathfrak{E})$ such that $\{ J_i
|_N \}_i$ form a coordinate system of $N$. Indeed, by Conjecture
\ref{conj}, which we're assuming, the Frobenius manifolds $P$ and
$N$ are isomorphic via $\varphi$.

First, we rewrite the PDE (\ref{J-eq}) as
\begin{equation}\label{pde2}
z\partial _i
\partial _j J = (\partial _i\star\partial _j) J.
\end{equation}
This is useful in computations.

Next, if $\xi$ and $\eta$ are $^\omega\nabla $ -
horizontal vector fields, then
\begin{eqnarray*} z\partial _{\xi\circ \eta}\partial
_\omega J_i &=& \partial _{(\xi\circ \eta)\star\omega} J_i \\
&=&
\partial _{\xi\star(\eta\star\omega)}J_i
\\&=& z\partial _{\xi}\partial _{\eta\star\omega }J_i \\ &=& z^2\partial
_{\xi}\partial _{\eta}
\partial _\omega J_i
\end{eqnarray*} since $\omega$ and $\eta\star\omega$ are
$\nabla$-horizontal. \end{proof}

\begin{Rmk}
Lemma \ref{d omega} reveals the relation between $\tilde{t}$ and
the $s=\varphi (t)$: \[ \tilde{t} = s + \sum _{n=0}^\infty
\frac{(-1)^{\epsilon ( \beta )}Q^{\beta}}{n!}\sum _{i, \
0\ne\tilde{\beta}\mapsto\beta}\gamma _i <\gamma ^i\cup\omega ,
\omega , \underbrace{s,...,s}_{n}>_{0,n+2,\tilde{\beta }},
\] where $\{\gamma ^j\cup\omega\}$ is the basis of $H^*(X//\T)^a$ dual to $\{\gamma _i\cup\omega\}$, that is,
$\int _{X//\T}\gamma _i \cup \omega \cup \gamma ^{j} \cup \omega  = \delta _i ^j$.
\end{Rmk}

Define, for $\tau \in N$,
\begin{equation}\label{Ifcn} \begin{split}
I(\tau ,z)&:= \left.\left((\prod _{\alpha \in \Phi_+} z\partial
_\alpha ) J_{X//T}\right)\right|_{Q^{\tilde{\beta}}=
(-1)^{\epsilon(\beta)}Q^\beta, N}(\tau ,z)\\&= \sum _{\beta}
(-1)^{\epsilon(\beta)}Q^\beta\sum_{\tilde{\beta}\mapsto\beta}
\prod _{\alpha \in \Phi_+} \left(c_1(L_\alpha
)+z\int_{\tilde{\beta}} c_1(L_\alpha )\right)
\left.J^{\tilde{\beta}}_{X//\T}\right|_{N}(\tau ,z)
\end{split}
\end{equation}
where
$\partial _\alpha$ is the ($\nabla$-flat) vector field associated to
$c_1(L_\alpha )$, the derivative of $J$ is taken
component-wise and $J^{\tilde{\beta}}_{X//\T}$ is the coefficient
of $Q^{\tilde{\beta}}$ in $J_{X//\T}$ before specializing the Novikov variables.
The latter equality follows from the divisor axiom.

\begin{Thm} \label{IJ}
There are unique $C^i (t,z)\in N(X//\G) [z][[t]]$ such that
$$I(\varphi (t),z) = \sum _{i}
C^i (t,z)z\partial _{t_i} \tilde{J}_{X//\G}(t,z)\cup \omega .$$
\end{Thm}

\begin{proof} For the proof we use Givental's
description \cite{SGF} of the rational Gromov-Witten theory for a
projective manifold $Y$ by means of a certain Lagrangian cone
$\mathcal{L}_Y$ with special properties (see Theorem 1 in \cite{SGF}).

Let
$s \in N$. By the very definition
$$I(s ,-z):=\pm\left(
(\prod _{\alpha \in \Phi_+} z\partial _\alpha)
J_{X//\T}\right)|_{Q^{\tilde{\beta}}=
(-1)^{\epsilon(\beta)}Q^\beta}(s,-z)\in zL ,$$ where $zL:=
zT_p\mathcal{L}_{X//\T}$ is the tangent space to the Lagrangian
cone at the point $p=J_{X//\T}(s)$.

Let $\{\phi _\mu\}$ be a basis of $H^*(X//\T)$ obtained by
adjoining to the basis $\{ \gamma _i \cup\omega \}$ of the
$\bW$-anti-invariant subspace $H^*(X//\T)^a$ a basis of
$(H^*(X//\T)^a )^\perp$. Since $\{z\partial _\mu
J_{X//\T}(s,-z)\}$ form a basis of $zL/z^2L$ over $N(X//\T)$,
$$I(s,z) =\sum C^\mu (s,z)z\partial _\mu J_{X//\T}|_{Q^{\tilde{\beta}}=
(-1)^{\epsilon(\beta)}Q^\beta}(s,z)$$ for some unique $C^\mu
(s,z)\in N(X//\T)[z][[s]]$.

\noindent Since $I$ is $\bW$-anti-invariant by construction,
the terms corresponding to the basis of
$(H^*(X//\T)^a )^\perp$ must vanish and we obtain
\begin{equation}
I(s,z )=\sum _i C^{\gamma _i \cup \omega } (s,z )(z
\partial _{\gamma _i\cup \omega }J_{X//\T})|_{Q^{\tilde{\beta}}=
(-1)^{\epsilon(\beta)}Q^\beta}(s,z).
\end{equation}
Now $\partial_{s_i}\star\omega=\partial_{\gamma_i\cup\omega}$,
therefore by equation (\ref{pde2})
\begin{multline}
\sum _i C^{\gamma _i \cup \omega } (s,z )(z
\partial _{\gamma _i\cup \omega }J_{X//\T})|_{Q^{\tilde{\beta}}=
(-1)^{\epsilon(\beta)}Q^\beta}(s,z)=\\
\sum _i C^{\gamma _i\cup \omega }(s,z) (z
\partial _{s_i}z\partial_{\omega }J_{X//\T})|_{Q^{\tilde{\beta}}=
(-1)^{\epsilon(\beta)}Q^\beta}(s,z).
\end{multline}
Finally, Lemma 5.3.1 gives
\begin{multline}
\sum _i C^{\gamma _i\cup \omega }(s,z) (z\partial
_{s_i}z\partial_{\omega }J_{X//\T})|_{Q^{\tilde{\beta}}
 =(-1)^{\epsilon(\beta)}Q^\beta}(s,z)=\\ \sum _i C^{\gamma _i \cup
\omega }(t,z) z
\partial _{t_i }\tilde{J}_{X//\G}(t,z)\cup \omega,
\end{multline}
where $\varphi (t)=s$. The Theorem follows from (5.3.3), (5.3.4)
and (5.3.5). \end{proof}

\begin{Cor} \label{J-I correspondence}
$$\tilde{J}_{X//\G} (t,z) \cup \omega = I(\tilde{t},z)+\sum _i C^i
(\tilde{t},z )z\partial _{\tilde{t} _i} I (\tilde{t} ,z)$$ for
some unique $C^i (\tilde{t},z)\in N(X//G) [[z,\tilde{t}]]$, where
$\tilde{t} = \sum \tilde{t} _i\gamma _i$. The expression of $t$ in
terms of $\tilde{t}$ is uniquely determined by the expansion of
the right-hand side as $z+t(\tilde{t} )+ O(z^{-1})$ (and coincides
with the formula (\ref{mirrortransf})).
\end{Cor}

\begin{proof} The theorem above shows that, with the identification of
cohomology spaces $H^*(X//\T)^{a}$ with $H^*(X//\G)$ by the map
$\tilde{\sigma} \cup \omega \mapsto \sigma$, the $I$-function
generates the Lagrangian cone $\mathcal{L}_{X//\G}$ describing the
rational Gromov-Witten theory of $X//\G$ \cite{SGF}. Since
$\{z\partial _{\tilde{t} _i} I (\tilde{t},-z )\}$ also form a
basis of $L/zL$, where $L$ is the tangent space of
$\mathcal{L}_{X//\G}$ at the point $I(\tilde{t},-z )$, the
Corollary follows.

A constructive argument may also be given, using the \lq\lq
Birkhoff factorization" method. See \cite{CG}, Corollary 5 and the
paragraph before it for details. \end{proof}

Corollary \ref{J-I correspondence} is a generalization of
Conjecture 4.3 in \cite{BCK2} to the \lq\lq big" parameter space.
The arguments in this section can be reversed to show that the
Corollary implies Conjecture \ref{conj}

\section{Flag manifolds for other classical types}\label{others}

In this section, we extend the
abelian/nonabelian correspondence in the
presence of additional twists by homogeneous vector bundles and
apply it to the case of generalized flag manifolds of Lie groups of types $B$,
$C$, $D$.

\subsection{Twisting by bundles}

Let $\bf S \times \G$ act on $X$ as in section \ref{Setting}. Let
$\mathcal{V}$ be a $\G$-representation space (as in [BCK2]).
There are $\bf S\times \G \times \CC ^*$-actions on $X$ and
$\cV$ (where $\CC ^*$ acts trivially on $X$ and homothetically on
$\cV$), inducing $\bf S\times \CC ^*$-equivariant vector bundles
$$\mathcal{V}_\T :=X^s({\bf T})\times _\T \cV , \ \ \
\mathcal{V}_\G:=X^s(\G )\times _\G \cV$$ over nonsingular
quotients $X//\T$ and $X//\G$, respectively. Put
$\CC[\lambda ']:=H^*(B\CC^*)$.

There is an $\bS \times \CC ^*$-equivariant
Frobenius structure on $$Z':=Spf ((N(X//\T)[\lambda] \otimes
_{\CC}\CC ((\frac{1}{\lambda '})) [[\left( H^*_\bS(X//\T)\otimes
(N(X//\T)[\lambda] \otimes _{\CC}\CC ((\frac{1}{\lambda
'})))\right)^\vee]]$$ defined by the $\bS\times \CC ^*$-equivariant
genus zero Gromov-Witten invariants
of $X//\T$ {\it twisted by (the equivariant Euler class of)} $\mathcal{V}_\T$. Here we
introduce the extra coefficient ring $\CC ((\frac{1}{\lambda '}))$
to invert
$$c_{\mathrm{top}}^{\bS\times\CC ^*}(\mathcal{V}_\T)=\sum
_{i=0}^{\mathrm{rk}\cV _\T} (\lambda ')^{\mathrm{rk}\cV _\T
-i}c_i^S (\cV _\T).$$
We list some comments on this Frobenius structure for
clarification, and refer the reader to \cite{CG} for details.

\begin{itemize}

\item The twisted metric $g_{\mathcal{V}_\T}$ is given by
$$g_{\cV _\T}(a,b):=\int_{X//\T}a \cup b \cup c_{\mathrm{top}}^{\bS\times\CC
^*}(\mathcal{V}_\T), \ \ \text{for } a,b\in H^*_{\bS}(X//\T).$$

\item The twisted product is given by the requirement that
\begin{eqnarray*} & & g_{\cV _\T}(a*_{\cV _\T}b,c)   =
\langle\langle a,b,c \rangle\rangle _{\cV _\T} \\ &:=& \sum
_{\tilde{\beta}\in NE_1(X//\T)}\sum_{n} \frac{Q^{\tilde{\beta}}}{n!}
\int _{ [\overline{M}_{0,n+3}(X//\T ,\tilde{\beta} )
]^{vir}}ev_1^*(a)ev_2^*(b)ev_3 ^*(c) \\ & & \ \ \ \ \ \ \
ev_4^*(t)\; ...\; ev_{n+3}^*(t)\; c_{\mathrm{vir. top}}^{\bS\times\CC
^*}(R^\bullet\pi _* ev_{n+4}^*\mathcal{V}_\T),
\end{eqnarray*}
where $\pi$ denotes the projection $\overline{M}_{0,n+4}(X//\T ,\tilde{\beta} )
\ra \overline{M}_{0,n+3}(X//\T ,\tilde{\beta} )$
of moduli stacks of stable maps which forgets the last marked point.

\item The Euler vector field is $\mathfrak{E}_{\mathcal{V}_\T} =
\mathfrak{E}+\cE_{\bS} +\cE_{\CC ^*}-c_1 ^{\bS}(\cV _\T)$.

%\item Note that this Frobenius is qc type by divisor axiom.

\item
%Denote by $J_{\cV _\T}^{\bS\times \CC ^*}$
The normalized ($\bS\times \CC ^*$-equivariant) $J$-function is
$$J_{\cV _\T}^{\bS\times \CC ^*}: t \mapsto z+t + \sum _i \phi ^i
\langle\langle\frac{\phi _i}{z-\psi}\rangle\rangle _{\cV _\T}, $$
where $\{\phi _i\}$ and $\{\phi ^i\}$ are dual bases with respect
to the twisted metric $g_{\cV _\T}$.

\end{itemize}

Similarly, we construct an $\bS\times \CC ^*$-equivariant
Frobenius structure on the formal scheme $P'$ associated to
$H^*_\bS(X//\G)\otimes (N(X//\G)[\lambda]\otimes \CC
((\frac{1}{\lambda '})))$ using genus zero $\bS\times \CC
^*$-equivariant Gromov-Witten invariants on $X//\G$ twisted by
$\mathcal{V}_\G$.

\bigskip

Now, as in section \ref{induced}, we can further twist the
Frobenius structure on $Z'$ by $\omega :=
\sqrt{\frac{1}{|\bW|}\prod _{\ka\in\Phi} c_1^\bS (L_\ka )} $ in
order to induce an $\bS\times \CC ^*$-equivariant Frobenius
structure on the formal scheme $N'$ over $N (X//\G)[\lambda]
\otimes _{\CC} \CC ((\frac{1}{\lambda '}))$ obtained as in loc.
cit. by fixing a lift of $H^*_{\bS}(X//\G )$ to $H^*_{\bS}(X//\T
)^{\bf W}$.

\begin{Conj}\label{conj twist} Let $\varphi: P' \ra N'$ be the isomorphisms of formal
schemes over $N(X//G)[\lambda ]\ot _{\CC}\CC ((\frac{1}{\lambda
'}))$ defined by $\varphi ^*(s_i)=t_i$. Then $\varphi $ induces an
isomorphism of formal $\bS\times\CC ^*$-equivariant Frobenius
structures.
\end{Conj}

\begin{Thm}\label{twist J} Conjecture \ref{conj} implies Conjecture \ref{conj
twist}, and furthermore, $$\tilde{J}_{\cV _\G}^{\bS\times \CC ^*}
(t,z) \cup\omega = z\partial _\omega J_{\cV _\T}^{\bS\times \CC
^*}|_{Q^{\tilde{\beta}} = (-1)^{\epsilon (\beta)} Q^\beta , N' }
(\varphi (t),z).$$
\end{Thm}

\begin{proof} It is enough to show the equality of $J$-functions
above, since it implies that $\varphi$ preserves the product
structures.

Abusing notation, for $\gamma \in
H^*_{\bS\times\CC ^*}(X//\T)^{a}$, $\sigma \in H^*_{\bS\times \CC
^*}(X//\G)$, denote $\sigma$ by $\frac{\gamma}{\omega}$  if
$\tilde{\sigma}\cup\omega =\gamma$. We also denote by $\mathcal{L}_{X//\G}^{\bS}$,
$\mathcal{L}_{X//\G}^{\bS\ti \CC ^*}$, and $\mathcal{L}_{\cV _\G}^{\bS\times \CC ^*}$
the Lagrangian cones given respectively by the $\bS$-equivariant,
$\bS\times\CC^*$-equivariant, and $\cV_\G$-twisted, $\bS\times\CC^*$-equivariant
rational GW invariants of $X//\G$.

By (the $\bS$-equivariant
version of) Lemma \ref{d omega}
$$\frac{z\partial _\omega J_{X//\T}^{\bS}|_{Q^{\tilde{\beta}}
= (-1)^{\epsilon (\tilde{\beta})} Q^\beta ,N' }(-z)}{\omega} \in
\mathcal{L}_{X//\G}^{\bS}.$$ Hence
$$\Delta _{\cV _\G} \frac{z\partial _\omega J_{X//\T}^{\bS}|_{
Q^{\tilde{\beta}} = (-1)^{\epsilon (\tilde{\beta})} Q^\beta , N'}
(-z) }{\omega} \in \Delta _{\cV _\G} \mathcal{L}_{X//G}^{\bS\times
\CC ^* } = \mathcal{L}_{\cV _\G}^{\bS\times \CC ^*}$$ by Corollary
4 in [CG], where
\begin{multline} \notag \Delta _{\cV _\G}=\prod _{\rho
_i :{\text{ Chern roots of }}\cV _G} b_{\rho _i} (\lambda ', z), \\
b_{\rho }(\lambda ', z)  = \exp \left(\frac{(\lambda '+\rho
)\ln (\lambda ' + \rho ) - (\lambda ' +\rho )}{z} + \sum _{m > 0}
\frac{B_{2m}}{2m (2m-1)}(\frac{z}{\lambda ' +\rho
})^{2m-1}\right)
\end{multline}
(and $B_{2m}$ are the Bernoulli numbers).
Since $$\Delta _{\cV _\G} \frac{z\partial _\omega
J_{X//\T}^{\bS}|_{ Q^{\tilde{\beta}} = (-1)^{\epsilon (\beta)}
Q^\beta , N'} }{\omega} = \frac{z\partial _\omega
\widetilde{\Delta}_{\cV _\G} J_{X//T}^{\bS}|_{Q^{\tilde{\beta}} =
(-1)^{\epsilon (\beta)} Q^\beta , N' }}{\omega}$$ and
$$\widetilde{\Delta}_{\cV _\G} = \Delta _{\cV _\T}
\ \ \ (\mathrm{mod}\ \mathrm{ker}(\cup \omega))$$ we conclude that
\begin{equation}\label{J_VG}\frac{z\partial _\omega J_{\cV _\T
}^{\bS\ti \CC ^*}|_{ Q^{\tilde{\beta}} = (-1)^{\epsilon (\beta)}
Q^\beta ,N'}(-z) }{\omega} \in \mathcal{L}_{\cV _\G}^{\bS\ti \CC
^*} .\end{equation}

Since the $J$-function $J_{\cV _\G} (-z)$ is uniquely
characterized by the intersection of the Lagrangian cone
$\mathcal{L}_{\cV _\G}$ with the subspace $-z+z\mathcal{H}_-$ as
in [Giv3], it follows that (\ref{J_VG}) is the $J$-function for
$P'$. That is,
$$J_{\cV _\G}^{\bS\ti \CC ^*} (t,z)
=\frac{z\partial _\omega J_{\cV _\T }^{\bS\ti \CC
^*}|_{Q^{\tilde{\beta}} = (-1)^{\epsilon (\beta)} Q^\beta
,N'}}{\omega} (\tau (t),z)$$ for some unique $\tau (t)$. As in
Corollary \ref{J-I correspondence}, the relation between $\tau$
and $t$ is given by the expansion of the right-hand side with
respect to $z$.

We have \begin{eqnarray*} g_{\cV _\G}(\partial _{t_i},\partial
_{t_j})+o(z) &=& g_{\cV _\G}(\partial _{t_i}J_{\cV _\G} ,\partial
_{t_j}J_{\cV _\G}) \\ & = & g_{\cV _\T}(z\partial _{t_i}\partial
_\omega J _{\cV _\T},z\partial _{t_i}\partial _\omega J_{\cV _\T})
\\
&=& g_{\cV _\T}(\partial _{\eta _i \star _{\cV
_\T}\omega},\partial _{\eta _j\star _{\cV _\T} \omega})+o(z)
\end{eqnarray*} where $\eta _i:=\partial _{t_i} (\tau ) $. We conclude that
$\eta _i \star _{\cV _\T}\omega = \gamma _i \cup \omega$, hence
$\tau (t)$ coincides with the map $\varphi$. \end{proof}

\begin{Rmk}\label{J VG} If $\cV _\G$ and $\cV _\T$ are
generated by $\bS$-equivariant global sections, then $J_{\cV
_\G}^{\bS}$ and $J_{\cV _\T}^{\bS}$ are well-defined without the
auxiliary variable $\lambda '$ (see \cite{CG}) and hence the
equality of $J$-functions in Theorem \ref{twist J} also holds
without $\lambda '$.
\end{Rmk}

\subsection{A simple Lemma}

Let $X$ be a nonsingular projective variety with an $\bS$-action
whose fixed points are isolated, and let $Y$ be a connected
component of the nonsingular zero locus of a regular $\bS$-equivariant
section of a $\bS$-equivariant bundle $E$. Suppose that $E$ is
generated by $\bS$-equivariant global sections. Let $i$ denote the
inclusion of $Y$ in $X$.

\begin{Lemma} If $i^*(\tilde{t})=t$, then $J_Y ^{\bS}(t,z)|_{Q^{\bf d}=Q^{i_*{\bf d}}}
= i^*J_E^{\bS} (\tilde{t},z)$ where $|_{Q^{\bf d}=Q^{i_*{\bf d}}}$
denotes the Novikov ring base change given by the pushforward
$i_*: NE_1(Y) \ra NE_1(X)$.
\end{Lemma}

\begin{proof} For each fixed point $p_i$ of $X$ under the $\bS$-action,
choose a nonzero class $\delta _i$ in $H^*_{\bS}(X)\otimes \CC
(\lambda )$ supported near $p_i$, and let $\{ \delta ^j \}$ be the
dual basis, that is, $\int _{c_{\mathrm{top}}^{\bS}(E)\cap [X]
}\delta ^i \cup \delta _j = \delta _{ij}$. Note that for nonzero
$\beta \in NE_1(X)$,
\begin{eqnarray*} i^*J_E  ^{\bS ,\beta}(\tilde{t},z) &=&
\sum _{k \ :\ p_k\in X^{\bS}} \frac{i^*\delta _k}{n!} \int
_{c_{\mathrm{top}}^{\bS}(\pi _*ev_{n+2}^*E)\cap
[\overline{M}_{0,n+1}(X,\beta)]^{vir}}
\frac{\delta ^k}{z-\psi} \prod _{i=1}^n ev_{1+i}^*(\tilde{t})  \\
&=& \sum _{k \ :\ p_k\in Y ^{\bS}}\frac{i^*\delta _k}{n!}\sum
_{{\bf d}\in NE_1(Y) \ : \ i_*{\bf d}=\beta} \int
_{[\overline{M}_{0,n+1}(Y,{\bf d} )]^{vir}}\frac{i^*\delta
^k}{z-\psi} \prod _{i=1}^n ev_{1+i}^*(t),
\end{eqnarray*} where the latter equality follows from \cite{KKP}.
Note that $i^*J_E  ^{\bS ,\beta}(\tilde{t})=0$ if there is no
${\bf d}\in NE_1(Y)$ such that $i_*{\bf d}=\beta$.
 Since $\{i^*\delta _k \}$
and $\{i^*\delta ^k \}$ form a dual pair of bases in
$H^*_{\bS}(Y)\otimes \CC (\lambda )$ with respect to the
equivariant Poincar\'e pairing, we are done. \end{proof}

\begin{Rmk}\label{J Y} The above Lemma is true for the nonequivariant
$J$-functions as well, since both sides of the identity can be
specialized to $\lambda =0$.
\end{Rmk}

\subsection{$J$-functions of flag manifolds of classical type}
Let $Y$ be a generalized flag manifold $K/P$,
with $K$ a simple complex Lie group of type $B$, $C$, or $D$
and $P$ a parabolic subgroup.
It can be viewed as a connected component of the zero
locus of a canonical section of a
homogeneous bundle $\cV _\G$ over an appropriate type $A$
partial flag manifold $X//\G = Fl(k_1,...,k_r,n)$. Here \[ \cV =\left\{
\begin{array}{ll} S^2(V^*) & \ \ \ \ \ \text{for types } B,\ D\\
\bigwedge ^2 V^*  & \ \ \ \ \ \text{for type } C
\end{array}\right. ,\]
where $V$ is the fundamental representation space of $GL_{k_r}(\CC
)$. Note that $\cV _\T$ is decomposable into a direct sum of line bundles
(since $\T$-representations are completely reducible).

Let $i:Y\subset X//\G$ be the natural inclusion and put
\begin{eqnarray*} I_{\cV _\G } &:=& \frac{1}{\omega}\left.\left((\prod _{\alpha
\in \Phi_+} z\partial _\alpha ) I_{\cV _\T
}\right)\right|_{Q^{\tilde{\beta}}=(-1)^{\epsilon (\beta )}Q^\beta, N'},\\
I_{\cV _\T} & := & \sum _{\tilde{\beta}\in NE_1(X//\T)}\prod
_{k=1}^{\int _{\tilde{\beta}}\rho _i} \prod _{\rho _i : \text{
Chern roots of } \cV _\T }(\rho _i + kz)
J_{X//\T}^{\tilde{\beta}}.
\end{eqnarray*}
Note that $I_{\cV _\T}$ is a $H^*(X//\T)$ - valued series and
$I_{\cV _\G}$ is a $H^*(X//\G)$ - valued series.

Let $\bS$ be a maximal abelian subgroup of the simple complex Lie group
$K$. It acts on the flag manifold
$Fl(k_1,...,k_r,n)$ with isolated fixed points and
$Y$ is an $\bS$-invariant submanifold. Since bundles $\cV_\G$ and
$\cV_\T$ are generated by $\bS$-equivariant global sections and
$i^*: H^*_{\bS}(X//\G)\ra H^*_{\bS}(Y)$ (as well as $i^*: H^*(X//\G)\ra
H^*(Y)$) is surjective, we obtain the following

\begin{Cor} Fix a subspace $N_Y$ of $H^*(X//\T)^{\bf W}$ which is a lift
of $H^*(Y)$ under the composite surjection $i^*\circ(\pi ^*)^{-1}\circ j^*$.
The $J$-function of $Y$ can be expressed as
$$J_Y (t,z) |_{Q^{\bf d}=Q^{i_*{\bf d}}} = I_{\cV_\G}(\tau ,z)+\sum _k C ^k (\tau ,z)
 i^* (\frac{z\partial _{\tilde{t}_k}I_{\cV _\G} (\tau ,z)}{\omega})$$
for some unique $C^k(\tau ,z )\in N(X//\G)[[z,\tau]]$, where
$\tilde{t}_k$ are coordinates of $N_Y$.
\end{Cor}

\begin{proof}  Due to Remark \ref{J Y}, $J_Y = i^* J_{\cV _\G}$. Moreover,
by Remark \ref{J VG},
$$J_{\cV _\G}= \frac{z\partial _\omega}{\omega} J_{\cV _T}|_{
Q^{\tilde{\beta}} = (-1)^{\epsilon (\beta)} Q^\beta ,N'}.$$ Now
apply the quantum Lefschetz Theorem of Coates and Givental
\cite{CG} and use a similar argument to the one in the proof of
Theorem \ref{IJ} to conclude that $i^* (\frac{I_{\cV
_\G}(-z)}{\omega})$ generates the Lagrangian cone $\mathcal{L}_Y$.
\end{proof}

\begin{Rmk} This in particular reproves the result on small
$J$-function of flag manifolds of types $B$, $C$, $D$ in [BCK2]. No coordinate
change is necessary for the explicit description of this small
$J|_{t_{\text{small}}}$.
\end{Rmk}

\section {Appendix: Multi-point GW-invariants of Grassmannians}

Recall from \S\ref{LP reconstruction} the notation
$$I_{n,\beta}(\gamma_1,\dots,\gamma_n)=(-1)^{\epsilon(\beta)}\sum_
{\tilde{\beta}\mapsto\beta}\langle\gamma_1,\dots,\gamma_n\rangle_{0,n,\tilde{\beta}}^{X//\T}.$$
Theorem \ref{mainthm}, together with equation
(\ref{new product}) (or, better, the equation (\ref{partials equal})),
imply that Gromov-Witten invariants of a flag manifold can be written in terms
of invariants of the corresponding toric variety $X//\T$ by a formula of the form
\begin{equation*}
\langle\sigma_{i_1},\dots ,\sigma_{i_n} \rangle_{0,n,\beta}^{X//\G}
=I_{n,\beta}(\tilde{\sigma}_{i_1},\dots \tilde{\sigma}_{i_{n-2}},
\tilde{\sigma}_{i_{n-1}}\cup\omega,\tilde{\sigma}_{i_n}\cup\omega)+\mathrm{correction}
\end{equation*}
where \lq\lq correction" is an expression involving invariants
$I_{n',\beta'}(...,\tilde{\sigma}_a\cup\omega,\tilde{\sigma}_b\cup\omega)$
with $n'\leq n$ and $\beta'\leq\beta$. Without going into too many
details, this can be seen as follows. Using the double bracket
notation for derivatives of Gromov-Witten potentials mentioned in
\S 2.1, one writes (\ref{partials equal}) as
$$\langle\langle\sigma_i,\sigma_j\rangle\rangle_{X//\G}(s)=
\langle\langle\tilde{\sigma}_i,\tilde{\sigma}_j\rangle\rangle_{X//\T}(\tilde{t}(s)),
$$
with $\tilde{t}(s)$ the inverse of the change of variables
(\ref{mirrortransf}). This is an equality of power series in
$s$-variables, and the formula for GW-invariants is obtained by
identifying the coefficients of monomials in the $s_j$'s. The
coefficient of an $s$-monomial in the power series
$\tilde{t}_k(s)$ can be explicitly expressed using the Lagrange
Inversion Formula (see \cite{GJ}, Theorem 1.2.9) in terms of the
coefficients of $\tilde{t}$-monomials of {\it lower} total degree
in the power series $s(\tilde{t})$ from (\ref{mirrortransf}).

The above discussion shows that the correction term will in general be quite complicated.
Moreover, while it is possible in principle to give an exact
expression, this will require the use of Lagrange inversion for
computing the inverse $\tilde{t}(s)$ of the coordinate change (\ref{mirrortransf}),
or, equivalently, the inverse (expressed in $s$-variables) of
the matrix of quantum multiplication with $\omega$ on a lift of $H^*(X//\G)$.

However, since flag manifolds are Fano of
index $\geq 2$, a different approach that uses Lemma \ref{small space}$(i)$
will allow us to reduce to
computing only the inverse of the matrix of {\it small} quantum
multiplication with $\omega$. In the case of Grassmannians, when the associated
abelian quotient is a product of projective spaces, it is an easy observation
that the small quantum product with $\omega$ is trivial (\cite{BCK1}, Lemma 2.4),
hence no matrix inversion is necessary. We present the derivation of closed formulae
for Grassmannians in this appendix.

Let $Gr:=Grass(k,n)$ be the Grassmannian of $k$-planes in $n$-space, thought
of as the GIT quotient $\Hom(\CC^k,\CC^n)//GL_k(\CC)$. The abelian quotient
is $\PP:=(\PP^{n-1})^k$. We consider the usual Schubert basis $\{\sigma_\lambda\}$
of $H^*(Gr,\CC)$, indexed by partitions $\lambda$ whose
Young diagrams fit in a $k\times(n-k)$ rectangle. We denote by $\mathcal{P}(k,n)$
the set of all such partitions. The intersection form in this basis is given by
$$\int_{Gr}\sigma_\lambda\cup\sigma_\mu=\delta_{\mu\lambda^{\vee}},$$
where $\lambda^{\vee}$ the complementary partition to $\lambda$ in
the $k\times(n-k)$ rectangle. The Grassmannian has Picard number
1, so the Novikov ring is $\CC[[Q]]$. On the other hand, the
Picard group of $\PP$ is isomorphic to $\ZZ^k$ and is generated by
$H_1,\dots,H_k$, with $H_j$ the pull-back of the hyperplane class
on the $j^{\mathrm{th}}$ factor. The Novikov ring of $\PP$ is
$\CC[[Q_1,\dots,Q_k]]$, and the specialization of Novikov
variables is $Q_i=(-1)^{k-1}Q$. In this case we also have a \lq\lq
canonical" lifting of a class on $Gr$ to a $\bW$-invariant class
on $\PP$ by taking
$$\tilde{\sigma}_\lambda=S_\lambda(H_1,\dots,H_k),$$
with $S_\lambda$ the Schur polynomial of the partition $\lambda$.
A curve class $\tilde{d}=(d_1,\dots,d_k)$ on $\PP$ is a lift of the curve
class $d$ on $Gr$ if and only if $\sum_{i=1}^kd_i=d$.
Finally, we have
$$\omega=\sqrt{\frac{(-1)^{\binom{k}{2}}}{k!}}\prod_{i<j}(H_i-H_j).$$

Let $\lambda^1,\dots,\lambda^l$ be (not
necessarily distinct) partitions. The generating
function for the $l$-point invariants of $Gr$ with $\sigma_{\lambda^i}$'s as insertions is
$$\left.\langle\langle\sigma_{\lambda^1},\dots,\sigma_{\lambda^l}\rangle\rangle_{Gr}
\right|_{t_{\mathrm{small}}}=
\sum_{d\geq 0}q^d\langle
\sigma_{\lambda^1},\dots,\sigma_{\lambda^l}\rangle_{0,l,d}^{Gr},$$
where $q^d=(Qe^{t_{\mathrm{small}}})^{d}$.
We start with three-point invariants. Let $\xi_\lambda$ be
the horizontal vector field (for the connection
$^\omega\nabla$) in $\Theta_N$ corresponding to
$\sigma_\lambda$ via the isomorphism $\varphi$
of Frobenius manifolds in Theorem \ref{mainthm}.
We have (cf. (\ref{new product}))
\begin{equation*}\begin{split}
\langle\langle\sigma_\lambda,\sigma_\mu,\sigma_\nu\rangle\rangle_{Gr}(t)
&=\xi_\lambda(\xi_\mu(\xi_\nu(F')))(\varphi(t))\\
&=
\langle\langle\hat{\xi}_{\lambda},\tilde{\sigma}_{\mu}\cup\omega,
\tilde{\sigma}_\nu\cup\omega\rangle\rangle_{\PP}|_{Q_{i}=(-1)^{k-1}Q,N}(\varphi(t))
\end{split}\end{equation*}
where $\hat{\xi}_\lambda$ is an extension of $\xi_\lambda$ to a vector field on $M$.
To unburden the notation,
this extension of vector fields will be understood when
necessary, and the same letter will be used for a vector field in $\Theta_N$,
or its extension to $\Theta_M$. Moreover, the specialization of Novikov variables
and the restriction to $N$ will be denoted by $\langle\langle\dots\rangle\rangle^{\bar{}}$.
Hence we rewrite the last
equation as
\begin{equation}\label{3-point}
\langle\langle\sigma_\lambda,\sigma_\mu,\sigma_\nu\rangle\rangle_{Gr}(t)
=\langle\langle\xi_{\lambda},\tilde{\sigma}_{\mu}\cup\omega,
\tilde{\sigma}_\nu\cup\omega\rangle\rangle_{\PP}^{\bar{}}(\varphi(t)).
\end{equation}
By Lemma \ref{small space} $(i)$ we get
$$\left.\langle\langle\sigma_\lambda,\sigma_\mu,\sigma_\nu\rangle\rangle_{Gr}(t)\right|_{t_{\mathrm{small}}}
=\left.\langle\langle\xi_{\lambda},\tilde{\sigma}_{\mu}\cup\omega,
\tilde{\sigma}_\nu\cup\omega\rangle\rangle_{\PP}^{\bar{}}\right|_{\tilde{t}_{\mathrm{small}}}.$$
From the relation $\xi_\lambda\star\omega
=\tilde{\sigma}_\lambda
\cup\omega$, and the fact that $\tilde{\sigma}_\lambda\star\omega
|_{\tilde{t}_{\mathrm{small}}}
=\tilde{\sigma}_\lambda\cup\omega$, we obtain
\begin{equation}\label{restriction2}
\left.\xi_\lambda\right|_{\tilde{t}_{\mathrm{small}}}=
\tilde{\sigma}_\lambda.\end{equation}
It follows that
\begin{equation*}\begin{split}
\langle\sigma_\lambda,\sigma_\mu,\sigma_\nu\rangle_{0,3,d}^{Gr}&=
I^{\PP}_{3,d}(\tilde{\sigma}_{\lambda},\tilde{\sigma}_{\mu}\cup\omega,
\tilde{\sigma}_\nu\cup\omega)\\
&=(-1)^{(k-1)d}\sum_{d_1+\dots+d_k=d}\langle
\tilde{\sigma}_{\lambda},\tilde{\sigma}_{\mu}\cup\omega,
\tilde{\sigma}_\nu\cup\omega\rangle^{\PP}_{0,3,(d_1,\dots,d_k)},
\end{split}\end{equation*}
an equation which was proved in \cite{BCK1}.

To obtain 4-point invariants we take the derivative of the relation (\ref{3-point})
and get
\begin{multline}\label{4-point}
\langle\langle\sigma_\pi,\sigma_\lambda,\sigma_\mu,\sigma_\nu\rangle\rangle_{Gr}(t)
=\xi_\pi(\langle\langle\xi_{\lambda},\tilde{\sigma}_{\mu}\cup\omega,
\tilde{\sigma}_\nu\cup\omega\rangle\rangle_{\PP}^{\bar{}}(\varphi(t)))\\
=\langle\langle\xi_\pi,\xi_{\lambda},\tilde{\sigma}_{\mu}\cup\omega,
\tilde{\sigma}_\nu\cup\omega\rangle\rangle_{\PP}^{\bar{}}(\varphi(t))+
\langle\langle\nabla_{\xi_\pi}\xi_{\lambda},\tilde{\sigma}_{\mu}\cup\omega,
\tilde{\sigma}_\nu\cup\omega\rangle\rangle_{\PP}^{\bar{}}(\varphi(t)),
\end{multline}
where $\nabla=\nabla^{\PP}$ is the connection on $M$. Since
$$0= {^\omega\nabla}_{\xi_\pi}\xi_\lambda\star\omega=
\nabla_{\xi_\pi}(\xi_\lambda\star\omega)=
(\nabla_{\xi_\pi}\xi_\lambda)\star\omega+
\sum_{a\in\mathcal{P}(k,n)}
\langle\langle\xi_\pi,\xi_\lambda,\omega,\tilde{\sigma}_a\cup\omega\rangle\rangle_{\PP}^{\bar{}}
(\tilde{\sigma}_{a^\vee}\cup\omega),
$$
it follows that
\begin{equation}\label{derivative}
\nabla_{\xi_\pi}\xi_\lambda=-\sum_{a\in\mathcal{P}(k,n)}
\langle\langle\xi_\pi,\xi_\lambda,\omega,\tilde{\sigma}_a\cup\omega\rangle\rangle _{\PP}^{\bar{}}
\xi_{a^\vee}.
\end{equation}
Combining with (\ref{4-point}) we find
\begin{multline}\label{4-point new}
\langle\langle\sigma_\pi,\sigma_\lambda,\sigma_\mu,\sigma_\nu\rangle\rangle_{Gr}(t)
=\langle\langle\xi_\pi,\xi_{\lambda},\tilde{\sigma}_{\mu}\cup\omega,
\tilde{\sigma}_\nu\cup\omega\rangle\rangle_{\PP}^{\bar{}}(\varphi(t))\\-
\sum_{a\in\mathcal{P}(k,n)}
\langle\langle\xi_\pi,\xi_\lambda,\omega,\tilde{\sigma}_a\cup\omega\rangle\rangle _{\PP}^{\bar{}}
\langle\langle\xi_{a^\vee},\tilde{\sigma}_{\mu}\cup\omega,
\tilde{\sigma}_\nu\cup\omega\rangle\rangle_{\PP}^{\bar{}}(\varphi(t)).
\end{multline}
Now we restrict to $t_{\mathrm{small}}$, using (\ref{restriction2}), to get
\begin{multline}\notag
\langle\sigma_\pi,\sigma_\lambda,\sigma_\mu,\sigma_\nu\rangle_{0,4,d}^{Gr}=
I^{\PP}_{4,d}(\tilde{\sigma}_\pi,\tilde{\sigma}_{\lambda},\tilde{\sigma}_{\mu}\cup\omega,
\tilde{\sigma}_\nu\cup\omega)\\
-\sum_{a\in\mathcal{P}(k,n)}\sum_{e+f=d}
I^{\PP}_{4,e}(\tilde{\sigma}_\pi,\tilde{\sigma}_{\lambda},\omega,
\tilde{\sigma}_a\cup\omega)
I^{\PP}_{3,f}(\tilde{\sigma}_{a^\vee},\tilde{\sigma}_{\mu}\cup\omega,
\tilde{\sigma}_\nu\cup\omega).
\end{multline}
The following remark is in order: while the left-hand side of the last formula is
manifestly invariant under permutations of the indices $\pi,\lambda,\mu$, and $\nu$, it
is not at all obvious that the right-hand side has this property. The invariance can, however,
be checked directly
using the splitting axiom for Gromov-Witten invariants, the vanishing
result in Lemma \ref{vanishing}, and the triviality of the small quantum
product with $\omega$.

Taking another derivative in (\ref{4-point new}) we get
\begin{equation*}
\begin{split}
\langle\langle\sigma_\rho,&\sigma_\pi,\sigma_\lambda,\sigma_\mu,\sigma_\nu\rangle\rangle_{Gr}(t)
= \langle\langle\xi_\rho,\xi_\pi,\xi_{\lambda},\tilde{\sigma}_{\mu}\cup\omega,
\tilde{\sigma}_\nu\cup\omega\rangle\rangle_{\PP}^{\bar{}}(\varphi(t))\\
+\langle\langle&\nabla_{\xi_\rho}\xi_\pi,\xi_{\lambda},\tilde{\sigma}_{\mu}\cup\omega,
\tilde{\sigma}_\nu\cup\omega\rangle\rangle_{\PP}^{\bar{}}(\varphi(t))
+\langle\langle\xi_\pi,\nabla_{\xi_\rho}\xi_{\lambda},\tilde{\sigma}_{\mu}\cup\omega,
\tilde{\sigma}_\nu\cup\omega\rangle\rangle_{\PP}^{\bar{}}(\varphi(t))\\
-\sum_{a}&{\Big(}
\langle\langle\xi_\rho,\xi_\pi,\xi_\lambda,\omega,\tilde{\sigma}_a\cup\omega\rangle\rangle _{\PP}^{\bar{}}
\langle\langle\xi_{a^\vee},\tilde{\sigma}_{\mu}\cup\omega,
\tilde{\sigma}_\nu\cup\omega\rangle\rangle_{\PP}^{\bar{}}+\\
&\langle\langle\nabla_{\xi_\rho}\xi_\pi,\xi_\lambda,\omega,\tilde{\sigma}_a\cup\omega\rangle\rangle _{\PP}^{\bar{}}
\langle\langle\xi_{a^\vee},\tilde{\sigma}_{\mu}\cup\omega,
\tilde{\sigma}_\nu\cup\omega\rangle\rangle_{\PP}^{\bar{}}+\\
&\langle\langle\xi_\pi,\nabla_{\xi_\rho}\xi_\lambda,\omega,\tilde{\sigma}_a\cup\omega\rangle\rangle _{\PP}^{\bar{}}
\langle\langle\xi_{a^\vee},\tilde{\sigma}_{\mu}\cup\omega,
\tilde{\sigma}_\nu\cup\omega\rangle\rangle_{\PP}^{\bar{}}+\\
&\langle\langle\xi_\pi,\xi_\lambda,\omega,\tilde{\sigma}_a\cup\omega\rangle\rangle _{\PP}^{\bar{}}
\langle\langle\xi_\rho,\xi_{a^\vee},\tilde{\sigma}_{\mu}\cup\omega,
\tilde{\sigma}_\nu\cup\omega\rangle\rangle_{\PP}^{\bar{}}+\\
&\langle\langle\xi_\pi,\xi_\lambda,\omega,\tilde{\sigma}_a\cup\omega\rangle\rangle _{\PP}^{\bar{}}
\langle\langle\nabla_{\xi_\rho}\xi_{a^\vee},\tilde{\sigma}_{\mu}\cup\omega,
\tilde{\sigma}_\nu\cup\omega\rangle\rangle_{\PP}^{\bar{}}{\Big)}
(\varphi(t)).
\end{split}
\end{equation*}
As above, we use (\ref{derivative}) to replace
the $\nabla_{\xi_\bullet}\xi_\bullet$ insertions, then restrict to $t_{\mathrm{small}}$ and obtain
the following formula for $5$-point invariants:
\begin{multline}\notag
\langle\sigma_\rho,\sigma_\pi,\sigma_\lambda,\sigma_\mu,\sigma_\nu\rangle_{0,5,d}^{Gr}=
I^{\PP}_{5,d}(\tilde{\sigma}_\rho,\tilde{\sigma}_\pi,\tilde{\sigma}_{\lambda},
\tilde{\sigma}_{\mu}\cup\omega,
\tilde{\sigma}_\nu\cup\omega)\\
-\sum_a\sum_{e+f=d}
{\Big(}I^{\PP}_{5,e}(\tilde{\sigma}_\rho,\tilde{\sigma}_\pi,\tilde{\sigma}_{\lambda},\omega,
\tilde{\sigma}_a\cup\omega)
I^{\PP}_{3,f}(\tilde{\sigma}_{a^\vee},\tilde{\sigma}_{\mu}\cup\omega,
\tilde{\sigma}_\nu\cup\omega) \\
+I^{\PP}_{4,e}(\tilde{\sigma}_\rho,\tilde{\sigma}_{\pi},\omega,
\tilde{\sigma}_a\cup\omega)
I^{\PP}_{4,f}(\tilde{\sigma}_{a^\vee},\tilde{\sigma}_\lambda,\tilde{\sigma}_{\mu}\cup\omega,
\tilde{\sigma}_\nu\cup\omega)\\
+I^{\PP}_{4,e}(\tilde{\sigma}_\rho,\tilde{\sigma}_{\lambda},\omega,
\tilde{\sigma}_a\cup\omega)
I^{\PP}_{4,f}(\tilde{\sigma}_{a^\vee},\tilde{\sigma}_\pi,\tilde{\sigma}_{\mu}\cup\omega,
\tilde{\sigma}_\nu\cup\omega)\\
+I^{\PP}_{4,e}(\tilde{\sigma}_\pi,\tilde{\sigma}_{\lambda},\omega,
\tilde{\sigma}_a\cup\omega)
I^{\PP}_{4,f}(\tilde{\sigma}_{a^\vee},\tilde{\sigma}_\rho,\tilde{\sigma}_{\mu}\cup\omega,
\tilde{\sigma}_\nu\cup\omega){\Big)}\\
+\sum_{a,b}\sum_{e+f+h=d}
{\Big(} I^{\PP}_{4,e}(\tilde{\sigma}_\rho,\tilde{\sigma}_{\pi},\omega,
\tilde{\sigma}_b\cup\omega)I^{\PP}_{4,f}(\tilde{\sigma}_{b^\vee},\tilde{\sigma}_\lambda,\omega,
\tilde{\sigma}_a\cup\omega)I^{\PP}_{3,h}(\tilde{\sigma}_{a^\vee},\tilde{\sigma}_{\mu}\cup\omega,
\tilde{\sigma}_\nu\cup\omega)\\
+I^{\PP}_{4,e}(\tilde{\sigma}_\rho,\tilde{\sigma}_{\lambda},\omega,
\tilde{\sigma}_b\cup\omega)I^{\PP}_{4,f}(\tilde{\sigma}_{b^\vee},\tilde{\sigma}_\pi,\omega,
\tilde{\sigma}_a\cup\omega)I^{\PP}_{3,h}(\tilde{\sigma}_{a^\vee},\tilde{\sigma}_{\mu}\cup\omega,
\tilde{\sigma}_\nu\cup\omega)\\
 +I^{\PP}_{4,e}(\tilde{\sigma}_\pi,\tilde{\sigma}_{\lambda},\omega,
\tilde{\sigma}_b\cup\omega)I^{\PP}_{4,f}(\tilde{\sigma}_{b^\vee},\tilde{\sigma}_\rho,\omega,
\tilde{\sigma}_a\cup\omega)I^{\PP}_{3,h}(\tilde{\sigma}_{a^\vee},\tilde{\sigma}_{\mu}\cup\omega,
\tilde{\sigma}_\nu\cup\omega){\Big)}
\end{multline}

It is now clear how to proceed to obtain and prove by induction a formula for
Gromov-Witten invariants with an arbitrary
number of insertions. We leave this to the reader.

\end{document}